\documentclass[10pt]{amsart}
\pdfoutput=1
\usepackage{amsmath, amsthm, amssymb,slashed,stmaryrd}
\usepackage{ifpdf}
\usepackage[pdftex]{graphicx}
\usepackage{tikz}
\usetikzlibrary{matrix,arrows,calc}

\usepackage[pdftex,plainpages=false,hypertexnames=false,pdfpagelabels]{hyperref}
\long\def\todo#1{{\color{red} {#1}}}

 \theoremstyle{plain}
 \newtheorem{thm}{Theorem}[section]
 \newtheorem{cor}[thm]{Corollary}
 \newtheorem{lem}[thm]{Lemma}
  
 \newtheorem{prop}[thm]{Proposition}
 \theoremstyle{definition}
 \newtheorem{defn}[thm]{Definition}
 
  \newtheorem{convention}[thm]{Convention}

 \newtheorem{ex}[thm]{Example}
 \theoremstyle{remark}
 \newtheorem{rmk}[thm]{Remark}

\def\beq{\begin{eqnarray}}
\def\eeq{\end{eqnarray}}

\DeclareSymbolFont{bbold}{U}{bbold}{m}{n}
\DeclareSymbolFontAlphabet{\mathbbold}{bbold}

 \newcommand{\bp}{\begin{proof}[Proof]}
 \newcommand{\ep}{\end{proof}}

\DeclareMathOperator{\Vect}{{\sf Vect}}

\def\Pf{{{\rm det}^{1/2}}}
\def\part{f}

\def\Hodge{\mathcal{L}}

\def\Mfld{{\sf Mfld}}

\def\Sh{{\sf Sh}}
\def\TMfld{T\hbox{-{\sf Mfld}}}
\def\Mfld{{\sf Mfld}}
\def\X{\mathcal{X}}
\def\Y{\mathcal{Y}}

\def\Mell{\mathcal{M}_{\rm ell}}
\def\EE{\mathcal{E}}

\def\Field{\Phi}

\def\Ell{{\rm Ell}{}}

\def\pt{{\rm pt}}
\def\Map{{\sf Map}}
\def\Wit{{\rm Wit}}

\def\odd{{\rm odd}}
\def\ev{{\rm ev}}

\def\R{{\mathbb{R}}}

\def\id{{{\rm id}}}

\def\MF{{\rm MF}}

\def\cH{{\mathbb{H}}}
\def\HH{{\mathbb{H}}}
\def\K{{\rm {K}}}
\def\KO{{\rm {KO}}}

\def\H{{\rm {H}}}
\def\C{{\mathbb{C}}}
\def\CP{{\mathbb{CP}}}

\def\T{{\mathbb{T}}}
\def\Z{{\mathbb{Z}}}

\def\End{\mathop{\sf End}}

\def\bS{{\mathbb{S}}}

\def\SL{{\rm SL}}

\def\TMF{{\rm TMF}}

\newcommand\triplerightarrow{
\mathrel{\substack{\textstyle\rightarrow\\[-0.5ex]
                      \textstyle\rightarrow \\[-0.5ex]
                      \textstyle\rightarrow}}
}
\newcommand\quadrightarrow{
\mathrel{\substack{\textstyle\rightarrow\\[-0.5ex]
\textstyle\rightarrow\\[-0.5ex]
                      \textstyle\rightarrow \\[-0.5ex]
                      \textstyle\rightarrow}}
}

\begin{document}

\title{Supersymmetric localization, modularity and the Witten genus}

\author{Daniel Berwick-Evans}

\address{Department of Mathematics, University of Illinois at Urbana-Champaign}

\email{danbe@illinois.edu}

\date{\today}

\begin{abstract}
Equivariant localization techniques give a rigorous interpretation of the Witten genus as an integral over the double loop space. This provides a geometric explanation for its modularity properties. It also reveals an interplay between the geometry of double loop spaces and complex analytic elliptic cohomology. In particular, we identify a candidate target for the elliptic Bismut--Chern character. 
\end{abstract}

\maketitle 
\setcounter{tocdepth}{1}
\tableofcontents

\section{Introduction and statement of results}

Witten's argument for the modularity of elliptic genera appeals to the path integral~\cite[\S2]{witten_dirac}. The goal of this paper is to explain his argument in terms of the geometry of double loop spaces. 

The main tool is \emph{supersymmetric localization}, an infinite-dimensional version of the Duistermaat--Heckman (or equivariant localization) formula. It allows one to define putative infinite-dimensional integrals in terms of finite-dimensional ones. The mathematical literature on this subject starts with the work of Witten~\cite[\S4]{susymorse} and Atiyah~\cite{AtiyahCircular} relating the index theorem in $\K$-theory with the geometry of loop spaces. This paper connects elliptic cohomology with the double loop space generalization.

In our setup, supersymmetric localization happens in two steps. The first is a quasi-isomorphism between an equivariant de~Rham complex of the double loop space of a manifold and the de~Rham complex of the manifold itself. The equivariant structure depends on an elliptic curve. Varying the curve leads to a sheaf of commutative differential graded algebras on the moduli stack of elliptic curves that is a model for complex analytic elliptic cohomology. The second step in supersymmetric localization defines the path integral as the side of the fixed point formula that is well-defined in the infinite-dimensional setting. The resulting integration map is a morphism of sheaves on the moduli stack of elliptic curves. On global sections, this constructs a cocycle refinement of the complexification of the Ando--Hopkins--Rezk--Strickland string orientation of topological modular forms (TMF). In particular, we obtain a geometric construction of the Witten genus as a modular form. 

\subsection{Prelude: Complex analytic elliptic cohomology} 
For an elliptic curve over~$\C$, elliptic cohomology is ordinary cohomology with coefficients in $\C[\beta,\beta^{-1}]$ where $|\beta|=-2$. For a smooth manifold $X$, we define a sheaf $(\Ell(X),d)$ of commutative differential graded algebras (cdgas) on the complex analytic moduli stack $\Mell\simeq \HH/\SL_2(\Z)$ of elliptic curves. This sheaf is a model for the complex analytic elliptic cohomology of~$X$ in the following sense. 

\begin{prop}\label{prop:1}
The sheaf $(\Ell(X),d)$ of cdgas on $\Mell$ has the following properties:
\begin{enumerate}
\item For $U\subset \cH\to \Mell$ an open subset of the upper half plane in~$\C$ classifying a complex analytic family elliptic curves, the cdga~$(\Ell(X)(U),d)$ computes the (de~Rham) cohomology of $X$ with coefficients in~$\mathcal{O}(U)[\beta,\beta^{-1}]$ for $|\beta|=-2$. 
\item The derived global sections of $(\Ell(X),d)$ (i.e., its hypercohomology) is $\TMF(X)\otimes \C$. 
\end{enumerate}
\end{prop}


Proposition~\ref{prop:1} shows that $(\Ell(-),d)$ gives a complex analytic version of the sheaf of $E_\infty$-ring spectra on the moduli stack of (algebraic) elliptic curves constructed by Goerss, Hopkins, and Miller~\cite{Goerss,GoerssHopkins,HopkinsICM94} whose local sections are elliptic spectra and global sections are~TMF. Of course, working over~$\C$ vastly simplifies both the construction of the sheaf and the computation of derived global sections. Our primary goal is to connect these structures in homotopy theory with ones in physics. Historically, such a connection began with the Witten genus~\cite{witten_dirac,Segal_Elliptic}. 

We start by recalling how the Witten genus is situated relative to elliptic cohomology. For a manifold $X$ with spin structure, there is a characteristic class $[\frac{p_1}{2}(TX)]\in \H^4(X;\Z)$. If this class vanishes, then $X$ has a \emph{string structure}. Ando, Hopkins, Rezk, and Strickland constructed the string orientation of TMF~\cite{AHS,AHR} which determines a wrong-way map $\pi_!$ for manifolds with string structure that fits into the following square,
\beq
\begin{tikzpicture}[baseline=(basepoint)];
\node (A) at (0,0) {$\TMF^k(X)$};
\node (B) at (5,0) {$\phantom{AAAAAAAa}$};
\node (BB) at (6.2,0) {$\TMF^k(X)\otimes \C\simeq \H^k(X;\MF)$};
\node (C) at (0,-1.5) {$\TMF^{k-d}(\pt)$};
\node (D) at (5,-1.5) {$\MF^{k-d}.$}; 
\draw[->] (A) to node [above] {$\otimes \C$} (B);
\draw[->] (A) to node [left] {$\pi_!$} (C);
\draw[->] (C) to node [below] {$\otimes \C$} (D);
\draw[->] (B) to node [right] {$\int_X(-)[\Wit(X)]$} (D);
\path (0,-.75) coordinate (basepoint);
\end{tikzpicture}\nonumber
\eeq
Above, $\MF$ is the graded ring of weakly holomorphic modular forms and $[\Wit(X)]\in \TMF^0(X)\otimes\C$ is the \emph{Witten class}. The integral of $[\Wit(X)]$ over $X$ is the \emph{Witten genus}. More generally, this class measures the difference between the complexification of $\pi_!$ and integration in ordinary cohomology valued in~$\MF$. In other words, the Witten class is the Riemann--Roch factor associated with the string orientation of $\TMF$. From Proposition~\ref{prop:1},~$[\Wit(X)]$ comes from a derived global section of $(\Ell(X),d)$. We obtain a cocycle representative by endowing $X$ with some additional structure. Namely for an oriented manifold~$X$, a \emph{geometric rational string structure} is a choice of Riemannian metric on~$X$ and~$H\in \Omega^3(X)$ so that $dH=p_1(TX)$ where $p_1(TX)=-\frac{1}{8\pi^2}{\rm Tr}(R^2) \in \Omega^4(X)$ for $R\in \Omega^2(X;\End(TX))$ the curvature of the Levi-Civita connection. 


\begin{prop}\label{prop:2} Let $X$ be a manifold with a chosen geometric rational string structure. The Witten class of $X$ has a unique cocycle refinement as a derived global section of $(\Ell(X),d)$, where the choice of geometric rational string structure is derived descent data.\end{prop}

In a bit more detail, the usual formula for $\Wit(X)$ of an oriented Riemannian manifold involves the 2nd Eisenstein series, which is not a modular form; see~\eqref{eq:WC}. As a consequence, the naive cocycle refinement fails to descend to a global section of $(\Ell(X),d)$. However, the data of a rational string structure fits together with $\Wit(X)$ to give a cocycle in a double complex computing the hypercohomology of the sheaf~$(\Ell(X),d)$. This determines a cocycle refinement of $[\Wit(X)]$ as a \emph{derived} global section. 


\subsection{Results 1: A double loop space model for elliptic cohomology} \label{sec:resultsI}
Let $T=\R^2/\Z^2$ be the standard 2-dimensional torus,~$X$ a smooth manifold, and $\Map(T,X)$ the double loop space of~$X$ (defined as a sheaf on the site of $T$-manifolds, see~\S\ref{sec:appen}). Consider the $\Z/2$-graded complex $(\Omega^\bullet(\Map(T,X))^{T},Q_\tau)$, of $T$-invariant differential forms for the precomposition action of~$T$ on itself. The differential $Q_\tau$ depends on a choice of complex structure~$T\simeq \C/(\Z\oplus\tau\Z)=E_\tau$. Then let $Q_\tau=d-\iota_{{\omega^\sharp_\tau}}$ where $d$ is the de~Rham differential and $\iota_{\omega^\sharp_\tau}$ is contraction with a complex vector field on $\Map(T,X)$. This vector field comes from the infinitesimal precomposition action, using the vector field $\omega^\sharp_\tau$ on $E_\tau\simeq T$ that is dual to the invariant 1-form $\omega_\tau$ on~$E_\tau=\C/(\Z\oplus \tau\Z)$ descending from~$dz$ on~$\C$.

\begin{thm}\label{thm:0} The natural restriction map 
\beq
(\Omega^\bullet(\Map(T,X))^{E_\tau},Q_\tau)\to (\Omega^\bullet(X),d)\label{eq:qiso1}
\eeq
is a quasi-isomorphism of $\Z/2$-graded chain complexes.
\end{thm}

One surprising aspect of this result is that $(\Omega^\bullet(\Map(T,X))^{E_\tau},Q_\tau)$ has no a priori reason to have Mayer--Vietoris sequences for open covers of $X$: there are maps $T\to X$ that do not factor through a given cover. However, in light of the~\eqref{eq:qiso1} we find that $(\Omega^\bullet(\Map(T,X))^{E_\tau},Q_\tau)$ is indeed local in $X$ in this sense. 

We promote~\eqref{eq:qiso1} to a quasi-isomorphism of sheaves on the moduli stack of elliptic curves as follows. We define a sheaf of cdgas~$(\Ell(\Map(T,X)),Q)$ on~$\Mell$ whose fiber at a curve~$E_\tau$ is the 2-periodic version of the $\Z/2$-graded complex $(\Omega^\bullet(\Map(T,X))^{E_\tau},Q_\tau)$. Restriction to the constant maps $X\subset \Map(T,X)$ determines a map of sheaves of cdgas on~$\Mell$
\beq
(\Ell(\Map(T,X)),Q)\stackrel{{\rm res}}{\to} (\Ell(X),d), \label{eq:restr}
\eeq
which at a fixed $\tau\in \HH$ is the 2-periodic version of the $\Z/2$-graded map~\eqref{eq:qiso1}.

\begin{thm}\label{thm:1}
The restriction map~\eqref{eq:restr} induces an isomorphism on hypercohomology for any suborbifold~$U/\Gamma\subset \Mell$ with $U\subset \cH$ and $\Gamma<\SL_2(\Z)$. 
\end{thm}

\begin{cor}\label{cor}
The sheaf $(\Ell(\Map(T,X)),Q)$ of cdgas on~$\Mell$ gives a cocycle model for the complex analytic elliptic cohomology of~$X$, i.e., this sheaf satisfies features (1) and (2) in Proposition~\ref{prop:1}.
\end{cor}

In physical language, sections of $\Ell(\Map(T,X))$ are classical observables for the~$\mathcal{N}=(0,1)$ supersymmetric sigma model with target~$X$, and the differential $Q$ is the the supersymmetry operator which depends on a complex structure on $T$. The first half of supersymmetric localization is Theorem~\ref{thm:0}: the $Q$-cohomology of classical observables localizes onto the constant maps. For constructions over the moduli stack $\Mell$, the study of $Q$-cohomology in physics (explained in greater detail in~\S\ref{sec:sigma}) motivates us to work with sheaves on $\Mell$ \emph{up to quasi-isomorphism}. Passing to the families-version of localization in Theorem~\ref{thm:1}, $Q$-cohomology is a geometric model for complex analytic elliptic cohomology. This explains one way in which the derived geometry familiar in the study of elliptic cohomology arises naturally from considerations in physics. 

\subsection{Results 2: Integration over the double loop space and the Witten genus}\label{results2}

The inclusion $X\hookrightarrow \Map(T,X)$ as the constant maps has an infinite rank normal bundle. Considering this normal bundle parameterized by the choices of complex structure on $T\simeq E_\tau$ results in a family of operators~$-[Q,Q]+\beta R$ on~$\Mell$. Here $R$ is the curvature of the Levi-Civita connection on~$X$ and $-[Q,Q]$ is a vector field on $\Map(T,X)$ inherited from the $E_\tau$-action under the choice of complex structure $T\simeq E_\tau$ on $T$. Mimicking equivariant localization formulas, supersymmetric localization of the path integral considers the integration map (see~\S\ref{sec:physics}) 
\beq
&&\langle f\rangle_\tau:=\int_{X} \frac{f|_X}{\det^{1/2}_{\rm ren}(-[Q,Q]+\beta R)},\quad \langle-\rangle\colon \Ell^\bullet(\Map(T,X)) \stackrel{\rm res}{\to} \Ell^\bullet(X)\to \Ell^{\bullet-d}(\pt)\label{eq:llint}
\eeq
that restricts to constant maps and applies the above (finite-dimensional) integral. In~\eqref{eq:llint}, $d={\rm dim}(X)$ and $\det^{1/2}_{\rm ren}(-[Q,Q]+\beta R)$ is a renormalized square root of the determinant (i.e., a \emph{Pfaffian}) defined as an infinite product of finite-dimensional Pfaffians. 

\begin{thm} \label{thm:int} The renormalized Pfaffian $\det^{1/2}_{\rm ren}(-[Q,Q]+\beta R)$ is a conditionally convergent product. Any choice of ordering the product constructs a derived global section of~$(\Ell(X),d)$ if and only if~$X$ has a rational string structure. Derived descent data is equivalent to the choice of a rational string structure on~$X$. The resulting map~\eqref{eq:llint} is a cocycle refinement of the pushforward associated with the string orientation of TMF. In particular the integral of~1 is the Witten genus of~$X$
$$
\langle 1\rangle=\int 1\cdot \Wit(X) \in \Gamma(\Mell;\Ell^{-d}(\pt))
$$ 
as a modular form of weight $d/2={\rm dim}(X)/2$. 
\end{thm}

In physical language, the path integral for the $\mathcal{N}=(0,1)$ supersymmetric sigma model localizes onto the constant maps and is modular if and only if the string anomaly vanishes. A choice of string structure cancels the anomaly; this choice is additional data. Modular invariance of the path integral then implies modular invariance of the partition function~$\langle 1\rangle$, which is the Witten genus. This is Witten's original argument for modularity~\cite[\S2]{witten_dirac}.

\subsection{Towards an elliptic Bismut--Chern character}

The Bismut--Chern character is a map~\cite{Bismut3,Bismut1}
\beq
{\rm BCh}\colon \Vect^\nabla(X)\to (\Omega^\bullet(LX)^{S^1},Q), \label{eq:BisCh}
\eeq
sending a vector bundle with connection to an even, $Q$-closed element in $\Omega^\bullet(LX)^{S^1}$ where $Q=d-\iota_{\partial_t}$ for $\partial_t$ the vector field generating loop rotation. 
Post-composition of~\eqref{eq:BisCh} with restriction to the constant loops gives the usual differential form-valued Chern character, 
$$
{\rm Ch}\colon \Vect^\nabla(X)\to \Omega^\ev(X),\qquad (V,\nabla)\mapsto {\rm Ch}(V,\nabla)={\rm Tr}(\exp(\nabla^2))\in \Omega^\ev(X)
$$
so that~\eqref{eq:BisCh} is a lift of the Chern character to the loop space. The $\Z/2$-graded chain complex $(\Omega^\bullet(LX)^{S^1},Q)$ is quasi-isomorphic to the $\Z/2$-graded de~Rham complex of $X$ (see Theorem~\ref{thm:JonesPetrack}), so~\eqref{eq:BisCh} does not encode any homotopical information beyond the ordinary Chern character. However, it does encode geometry, e.g., traces of holonomies. Furthermore, it reveals an interplay between loop space geometry and the index theorem: essentially all the steps in the proof of the index theorem have a loop space interpretation (e.g., see~\cite{Bismut2}).

It is tempting to expect an \emph{elliptic} Bismut--Chern character 
\beq
&&\left\{\begin{array}{c} { \rm Geometric \ cocycles \ for} \\ {\rm elliptic \ cohomology \ of} \ X \end{array} \right\}\stackrel{?}{\longrightarrow} (\Ell(\Map(T,X)),Q)\label{eq:ellBC}
\eeq
from cocycles for the elliptic cohomology of $X$ (over $\Z$) to $Q$-closed sections of $\Ell(\Map(T,X))$. This temptation comes partly from the striking parallels in the corresponding stories from physics; see~\S\ref{sec:AS} and~\S\ref{sec:sigma} below. Indeed, partition functions of 2-dimensional supersymmetric field theories over~$X$ admit a map like~\eqref{eq:ellBC}, suggesting that such field theories could provide a geometric model for elliptic cohomology. This idea is an old one, going back to Segal~\cite{Segal_Elliptic} and Stolz--Teichner~\cite{ST04}. In his thesis, Fei Han gave an interpretation of the usual Bismut--Chern character in K-theory within the Stolz--Teichner framework~\cite{Han}; see also the announced result~\cite[Theorem 1.19]{ST11}. A study of the situation one dimension higher relating 2-dimensional field theories to the geometry of double loop spaces will help illuminate the ingredients that go into the long-sought index theorem for~TMF. 

\subsection*{Outline}

In~\S\ref{sec:physics} we overview the ideas from physics. This section is independent from the remainder of the paper and can be skipped by readers less interested in this motivation. We have included it to explain where our ideas came from, describe the relation between constructions in this paper with analogous structures in K-theory, and illuminate possible generalizations (e.g., proposed higher genus variants as in~\cite{AlvarezSinger}). 

The mathematical exposition begins in~\S\ref{sec:complexanal}, where we provide the basic definitions for complex analytic elliptic curves and elliptic cohomology. In~\S\ref{sec:double} we define the sheaf $(\Ell(\Map(T,X)),Q)$ and prove Theorems~\ref{thm:0}  and~\ref{thm:1}. In~\S\ref{sec:integral}, we construct the localized path integral and prove Theorem~\ref{thm:int}. In~\S\ref{sec:appen} we give the sheaf-theoretic point of view on ($T$-invariant) differential forms on~$\Map(T,X)$. 

\subsection*{Conventions}

Throughout, $X$ denotes a closed, oriented manifold. The compactness condition can be loosened, but it complicates the constructions. We take smooth functions and differential forms to be $\C$-valued. Unless stated otherwise, tensor products of Fr\'echet spaces (e.g., algebras of functions and sections of vector bundles) use the projective tensor product. This is a completion of the algebraic tensor product with the important properties
\beq
&&C^\infty(M\times N)\simeq C^\infty(M;C^\infty(N))\simeq C^\infty(M)\otimes C^\infty(N)\quad C^\infty(N;V)\simeq C^\infty(N)\otimes V\label{eq:projten}
\eeq
for manifolds $M$ and $N$ and a Fr\'echet vector space $V$. All sheaves are strict, i.e., non-homotopical. For example, a sheaf of cdgas is a chain complex of sheaves with a compatible graded commutative multiplication. 

\subsection*{Acknowledgements} 

I thank Matt Ando, Kevin Costello, Theo Johnson-Freyd, Charles Rezk, Nat Stapleton, Constantin Teleman, and Arnav Tripathy for helpful conversations.

\section{Motivation from physics: Path integrals and localization}\label{sec:physics}

Supersymmetric localization is inspired by the Duistermaat--Heckman formula~\cite{DuistermaatHeckman} that gives conditions under which the stationary phase approximation to an integral is exact. Stationary phase techniques are commonly used in quantum field theory to approximate path integrals. It is then of great interest to find quantum field theories in which the stationary phase approximation to the path integral is (formally) exact. 

This fortuitous situation often arises in the presence of supersymmetry, where the supersymmetry algebra allows one to formally apply the proof of the Duistermaat--Heckman formula from equivariant cohomology~\cite{BerlineVergneEquiv,AtiyahBott}. The literature on supersymmetric localization is vast; the recent volume~\cite{localizationvol} gives a comprehensive overview, with the chapter~\cite{localizationvolintro} being particularly relevant to the methods in this paper. 

In this section we overview the key ideas of supersymmetric localization, with a focus on two examples:  $\mathcal{N}=1$ supersymmetric mechanics (related to $\K$-theory and the index theorem) and the $\mathcal{N}=(0,1)$ supersymmetric sigma model (related to elliptic cohomology). Our discussion of the physics is intensionally informal; one can make the setup mathematically rigorous by formulating the space of fields below as a stack on the site of supermanifolds. However, this language of \emph{superstacks} is not very widely known and can be quite technical. The geometry of our mathematical applications is more easily understood in standard language (e.g., differential forms) after one has translated the supersymmetric physics into the geometry of double loop spaces. This section explains the translation.

 \subsection{Localization in finite dimensions}\label{sec:fd}
There are a variety of closely related results concerning (equivariant) localization, e.g.,~\cite{DuistermaatHeckman,susymorse,AtiyahBott,Bismut1}. We focus on results of Berline and Vergne~\cite{BerlineVergne,BerlineVergne2}. Let~$M$ be a compact manifold with the action by a torus~$\T$. Identify an element $\xi\in\mathfrak{t}$ in the Lie algebra of~$\T$ with its corresponding vector field on~$M$. Consider the complex $(\Omega(M)^\xi,Q)$ of $\xi$-invariant differential forms on~$M$ with the differential $Q=d-\iota_\xi$, where $\iota_\xi$ is the contraction operator. By the Cartan formula, $\frac{1}{2}Q^2=[Q,Q]=-\mathcal{L}_\xi$ is the Lie derivative, so $Q^2=0$ on $\Omega(M)^\xi$. Let $M_0\subset M$ be the submanifold consisting of the zeros of the vector field~$\xi$. 
 
 \begin{thm}[\cite{BerlineVergne2}, Proposition~2.1] \label{thm:BV1} The natural restriction map 
 $$
 (\Omega(M)^\xi,Q)\to (\Omega(M_0),d)
 $$
 is a quasi-isomorphism of $\Z/2$-graded complexes. 
 \end{thm}
 
Berline and Vergne prove a related integration formula. Choose a $\T$-invariant metric on~$M$. Let~$\nu$ be the $\T$-equivariant normal bundle for the inclusion $M_0\subset M$ and $R\in \Omega^2(M_0;\End(\nu))$ be the curvature 2-form of~$\nu$ from the Levi-Civita connection on~$M$.
 
 \begin{thm}[\cite{BerlineVergne2}, Theorem~2.8] \label{thm:BV2}  For $\alpha\in \Omega(M)^\xi$ with $Q\alpha=0$, we have the equality
\beq
\int_M \alpha=\int_{M_0} \frac{\alpha}{\det^{1/2}(\mathcal{L}_\xi+R)}\label{eq:BV2}
 \eeq
where $\mathcal{L}_\xi=-[Q,Q]$ is the Lie derivative acting on~$\nu$. 
 \end{thm}
 
When $M_0$ is discrete, the integral~\eqref{eq:BV2} becomes a finite sum. Hence, a potentially complicated analytical problem (the integral) collapses to a purely algebraic one (computing determinants). When confronted with the path integral in quantum field theory, we might ask for similarly fortuitous situations that reduce complicated integrals to algebra. 

\subsection{Localization in supersymmetric quantum field theory}\label{sec:local}
The basic input data of a classical field theory is a \emph{space of fields} $\Field$ and a \emph{classical action functional} $\mathcal{S}\in C^\infty(\Field)$. One looks to compute \emph{quantum expectation values}
\beq
\langle \alpha \rangle \stackrel{?}{=}\int_\Field \alpha e^{-\mathcal{S}}[d\phi], \qquad \langle -\rangle\colon C^\infty(\Field)\to \C\label{eq:pathintegral}
\eeq
for \emph{classical observables} $\alpha \in C^\infty(\Field)$, where~$[d\phi]$ is the notoriously ill-defined measure on~$\Field$. 
\begin{ex}\label{Ex:sigma}
The \emph{2-dimensional $\sigma$-model} has fields $\Field=\Map(\Sigma,X)$ for $\Sigma$ a Riemann surface and $X$ a Riemannian manifold. The classical action is the energy, $\mathcal{S}(\phi)=\frac{1}{2}\int_\Sigma \|d\phi\|^2$, for $\phi\in \Map(\Sigma,X)$. There is no satisfactory definition of~\eqref{eq:pathintegral} for this example. 
\end{ex}
Instead of focusing on the measure $[d\phi]$, a common approach is to construct a map $\langle -\rangle\colon C^\infty(\Field)\to \C$ satisfying \emph{algebraic} properties mimicking those expected from the would-be integral. For example, 
the \emph{Batalin--Vilkovisky formalism} codifies algebraic properties of an equivariant de~Rham complex of~$\Field$ incorporating symmetries of the theory assuming the ``top form" $[d\phi]$ exists~\cite{BV1,BV2}. Our approach is similar in spirit: we show that in good situations, natural algebraic constraints reduce the construction of the map~\eqref{eq:pathintegral} to a certain finite-dimensional integral.
 
These algebraic constraints come from symmetries. \emph{Infinitesimal symmetries} are a Lie algebra acting on~$C^\infty(\Field)$ preserving~$\mathcal{S}$. In a \emph{supersymmetric} field theory, $C^\infty(\Field)$ has a $\Z/2$-grading and infinitesimal symmetries form a Lie \emph{superalgebra}. This superalgebra frequently contains an odd element~$Q$ and an even element~$P$ satisfying~\cite[\S1]{susymorse} 
$$
[Q,Q]=\frac{1}{2}Q^2=P,\qquad [Q,P]=0.
$$
In this case, the algebra of $P$-invariant observables $C^\infty(\Field)^{P}$ is a $\Z/2$-graded chain complex with differential~$Q$. The hypothetical measure $[d\phi]$ in~\eqref{eq:pathintegral} is expected to be invariant under symmetries, which would imply $\langle Q f\rangle=0.$ Hence, a map~\eqref{eq:pathintegral} on the $Q$-closed observables only depends on the $Q$-cohomology of $(C^\infty(\Field)^{P},Q)$. Often a quasi-isomorphic complex comes from the restriction of functions to a finite-dimensional subspace~$\Field_0\subset \Field$. Then any path integral~\eqref{eq:pathintegral} on $Q$-closed observables that is compatible with symmetries can be constructed as a finite-dimensional integral over~$\Field_0$: the path integral \emph{localizes} to~$\Field_0$. It is common practice in the physics literature to \emph{define} the path integral in terms of localization, e.g., see~\cite[pg.~4]{localizationvolintro}. 


\subsection{Localization in supersymmetric quantum mechanics}\label{sec:AS}

In part of their study of the index theorem, Atiyah and Witten~\cite{AtiyahCircular,susymorse} observed that the stationary phase approximation to a certain loop space integral recovers the $\hat{A}$-genus. A physical explanation of this comes from supersymmetric quantum mechanics with $\mathcal{N}=1$ supersymmetry. Supersymmetric classical mechanics takes as input a Riemannian manifold~$X$, and studies the space $\Field_{\sigma,X}$ consisting of pairs $(x,\psi)$
\beq
x\colon S^1\to X,\quad \psi\in \Gamma(S^1,\bS\otimes x^*TX),\qquad (x,\psi)\in \Field_{\sigma,X}\label{Eq:dim1}
\eeq
where $\sigma$ is a spin structure on~$S^1$ with spinor bundle $\bS$. Equip $\Field_{\sigma,X}$ with the functional
$$
\mathcal{S}(x,\psi)=\frac{1}{2}\int_{S^1} (\langle \dot x,\dot x\rangle +\langle \psi,\nabla_{\partial_t} \psi\rangle )dt
$$
for $t$ the standard coordinate on $S^1=\R/\Z$. If~$\bS$ is a trivial line bundle (equivalently, $\sigma$ is the \emph{odd} spin structure), a choice of trivialization identifies~$\Field_{\sigma,X}$ with the tangent space to the free loop space. The $\Z/2$-graded physics involving supersymmetry requires we view sections of~$\bS$ as odd, meaning a choice of trivialization of $\bS^*$ identifies~$\Field_{\sigma,X}$ with the ``odd" tangent bundle of~$LX=\Map(S^1,X)$. The upshot of this parity convention is that functions on $\Field_{\sigma,X}$ are differential forms on the loop space,\footnote{We recall that functions on the vector space $V$ when regarded as an \emph{odd} vector space are the exterior algebra of the dual vector space, $\Lambda^\bullet V^*$, as a $\Z/2$-graded algebra, and $\Omega^\bullet(LX)=\Gamma(\Lambda^\bullet T^* LX)$.} 
\beq
&&C^\infty(\Field_{\sigma,X})\simeq \Omega^\bullet(LX)\qquad {\rm (isomorphism \ depends \ on \ a\ trivialization \ of \ }\bS^*{\rm )}.\label{eq:loopfield}
\eeq
A trivialization of~$\bS^*$ also determines a \emph{supersymmetry}, meaning an (odd) vector field $Q$ on $\Field_{\sigma,X}$ that preserves the action functional~$\mathcal{S}$. Under the identification above,~$Q=d-\iota_{\partial_t}$ is an odd derivation on $\Omega^\bullet(LX)$ where the vector field $\partial_t$ is the generator for the $S^1$-action on $LX$ associated with loop rotation and $\iota_{\partial_t}$ is contraction with $\partial_t$. Hereafter, we assume~$\bS$ is trivializable so that an isomorphism~\eqref{eq:loopfield} can be chosen. 

Quantizing this classical field theory in the path integral formalism looks to evaluate
\beq
\langle \alpha\rangle\stackrel{?}{=}\int_{\Field_{\sigma,X}} \alpha e^{-\mathcal{S}}[dxd\psi],\qquad \alpha\in C^\infty(\Field_{\sigma,X})\simeq \Omega^\bullet(LX).\label{eq:path1}
\eeq
Restricting to $Q$-closed functions and requiring $\langle Q\alpha\rangle=0$, any possible definition of~\eqref{eq:path1} only depends on the $Q$-cohomology of~$(\Omega^\bullet(LX)^{S^1},Q)$. 

\begin{thm}[\cite{JonesPetrack}, Theorem~2.1]\label{thm:JonesPetrack} The restriction to the constant loops 
$$
(\Omega^{\ev/\odd}(LX)^{S^1},Q)\to (\Omega^{\ev/\odd}(X),d)
$$
is a quasi-isomorphism of $\Z/2$-graded chain complexes. 
\end{thm}

Hence, the path integral~\eqref{eq:path1} localizes to the constant loops $X\subset LX$, and a choice of volume form on~$X$ determines one (formally) on~$LX$. To identify the ``correct" choice, we formally apply Theorem~\ref{thm:BV2}. The normal bundle~$\nu$ to the inclusion $X\subset LX$ of the constant loops is an infinite rank bundle. The $S^1$-action on~$LX$ endows the fibers of~$\nu$ with an $S^1$-action, leading to an operator $-[Q,Q]+R=\mathcal{L}_{\partial_t}+R$ that acts on sections of $\nu$ where $\mathcal{L}_{\partial_t}$ is the Lie derivative and $R$ is the curvature, viewed as an endomorphism-valued 2-form. There are various ways to make sense out of the regularized Pfaffian of this operator, e.g., \cite[Equation~4.6]{AtiyahCircular}, \cite[\S3.7]{Bismut2} or~\cite[Theorem~5.1]{JonesPetrack}. 

\begin{thm}\label{thm:A}
The regularized Pfaffian of $\mathcal{L}_{\partial_t}+R$ is a cocycle representative of the $\hat{A}$-class of $X$, i.e., $\det^{-1/2}_{\rm ren}(\mathcal{L}_{\partial_t}+R)=\hat{A}(X)\in \Omega^\ev(X)$. 
\end{thm}

One then defines the path integral over~$LX$ as
\beq
&&\langle\alpha\rangle=\int_{\Field_{\sigma,X}}\alpha e^{-\mathcal{S}}[dxd\psi]:=\int_X \alpha|_X\wedge \hat{A}(X),\qquad \alpha\in \Omega(LX)^{S^1}, Q\alpha=0.\label{eq:LXint}
\eeq
The right-hand side is the pushforward in $\KO\otimes \C$ inherited from the analytic index in real K-theory. In particular, $\langle 1\rangle$ is the $\hat{A}$-genus of $X$. 

\begin{rmk} One can also quantize supersymmetric mechanics in the \emph{Hamiltonian formalism} using a version of geometric quantization. This quantization is generally obstructed, but when $X$ is spin the resulting quantum theory encodes the Dirac operator acting on spinors on $X$. The agreement between this canonical quantization and the path integral~\eqref{eq:LXint} is the physics argument for the Atiyah--Singer index theorem~\cite{susymorse,Alvarez}. \end{rmk}

\subsection{Localization in 2-dimensional supersymmetric sigma models} \label{sec:sigma}
Fix a spin Riemann surface $\Sigma$ with spinor bundle $\bS$ and a Riemannian manifold $X$. \emph{The space of fields} for the $\mathcal{N}=(0,1)$ sigma model with source $\Sigma$ and target $X$, denoted $\Field_{\Sigma,X}$, consists of 
\beq
(x,\psi)\in \Field_{\Sigma,X},\qquad x\colon \Sigma\to X,\quad \psi\in \Gamma(\Sigma;\overline{\bS}\otimes x^*TX)\label{eq:fields}
\eeq
where $x$ is a smooth map and $\psi$ is a (anti-chiral) spinor valued in the pullback tangent bundle. We recall that $\overline{\bS}\otimes \overline{\bS}\simeq \Omega^{0,1}_\Sigma$~\cite{AtiyahSpin}. \emph{The classical action} is the function on fields
$$
\mathcal{S}_{\Sigma}(x,\psi)=\frac{1}{2}\int_\Sigma (\langle \partial x,\bar\partial x\rangle +\langle \psi,\partial_\nabla \psi\rangle)\qquad \mathcal{S}_{\Sigma}\in C^\infty(\Field_{\Sigma,X})
$$
where $\partial_\nabla$ is the $\partial$-operator on $\Sigma$ twisted by the pullback of the Levi-Civita connection $\nabla$. 

A choice of trivialization of $\overline{\bS}^*$ identifies~$\Field_{\Sigma,X}$ with the (odd) tangent bundle of $\Map(\Sigma,X)$, 
\beq
&&C^\infty(\Field_{\Sigma,X})\simeq \Omega^\bullet(\Map(\Sigma,X))\qquad {\rm (isomorphism\ depends\ on\ a\ trivializaton\ of \ } \overline{\bS}^*{\rm )}\label{eq:fieldsdoubleloop}
\eeq
analogous to~\eqref{eq:loopfield}. A trivialization of $\overline{\bS}^*$ exists if and only if $\Sigma\simeq \C/\Z\oplus \tau \Z$ is genus~1 and the spin structure is odd (alias: periodic-periodic). Hereafter we assume that $\Sigma$ is of this form so that an isomorphism~\eqref{eq:fieldsdoubleloop} can be chosen. The choice also defines a  supersymmetry, which under~\eqref{eq:fieldsdoubleloop} is the odd derivation~$Q_\tau=d-\iota_{\omega^\sharp_\tau}$, for $\omega^\sharp_\tau$ the complex vector field on $\C/\Z\oplus \tau \Z$ dual to the 1-form~$\omega_\tau$ that descends from~$dz$ on~$\C$. The precomposition action of the elliptic curve $E_\tau= \C/\Z\oplus \tau \Z$ on $\Sigma$ promotes the vector field $\omega^\sharp_\tau$ to one on~$\Map(\Sigma,X)$. 

We again look to evaluate the path integral, 
\beq
\langle \alpha \rangle_{\tau}\stackrel{?}{=}\int\alpha e^{-\mathcal{S}_{\Sigma}}dxd\psi\qquad \alpha\in C^\infty(\Field_{{\Sigma},X})\simeq \Omega^\bullet(\Map(\Sigma,X)).\label{eq:pathint}
\eeq
In parallel to Theorem~\ref{thm:A}, Theorem~\ref{thm:0} shows that the restriction map
$$
(\Omega^{\ev/\odd}(\Map(\Sigma,X))^{E_\tau},Q_\tau)\to (\Omega^{\ev/\odd}(X),d)
$$
is a quasi-isomorphism. Hence, we may define the localized path integral analogously to~\eqref{eq:LXint}. Theorem~\ref{thm:int} computes the path integral defined this way and compares it to the pushforward in elliptic cohomology associated to the string orientation. The agreement is analogous to how~\eqref{eq:LXint} yields the pushforward in complexified K-theory associated with the spin orientation. 

A new feature in the 2-dimensional case is that we can refine the constructions to ones over the moduli stack of elliptic curves. The isomorphism~\eqref{eq:fieldsdoubleloop} depends on a holomorphic trivialization of the dual spinor bundle~$\overline{\bS}^*\simeq \bS$. This is the data of a nonzero section of the square root of the Hodge bundle,~$\Hodge^{1/2}$, on $\Mell^{\rm Spin}$. Hence we have a \emph{canonical} isomorphism of $\Z/2$-graded algebras,
\beq
C^\infty(\Field_{\Sigma,X})\simeq \prod_{k\in \Z} \Hodge^{\otimes k/2}\otimes \Omega^k(\Map(\Sigma,X)).\label{eq:Z2}
\eeq
The differential $Q_\tau$ also depends on the trivialization of~$\bS$: $Q$~maps between different tensor powers of $\Hodge^{\otimes k/2}$. This prompts us to consider the $\Z$-graded version with $m$th graded piece
\beq
\Gamma(\Field_{\Sigma,X}; \Hodge^{\otimes -m/2})^{E_\tau} \simeq \Big(\prod_{k\in \Z} \Hodge^{\otimes (k-m)/2}\otimes \Omega^k(\Map(\Sigma,X))\Big)^{E_\tau},\label{eq:chainvect}
\eeq
where the differential takes the form~$Q_\tau=d-\omega_\tau^{-1}\iota_{\omega_\tau^\sharp}$, and depends on the data of the vector field $\omega^\sharp_\tau$ and the trivializing section~$\omega_\tau^{-1}$ of the fiber of~$\Hodge^{-1}=\Hodge^{\otimes -2/2}$ at $E_\tau$. Next we observe that the graded vector space $\Omega^\bullet(\Map(\Sigma,X))^{E_\tau}$ only depends on $\Sigma=\C/(\Z\oplus\tau\Z)$ as a smooth manifold, and a choice of diffeomorphism $\Sigma\simeq \R^2/\Z^2=T$ gives an isomorphism
\beq
\Gamma(\Field_{\Sigma,X}; \Hodge^{\otimes-\bullet/2})^{E_\tau} \simeq \prod_{k\in \Z} \Hodge^{\otimes (k-\bullet)/2}\otimes \Omega^{k}(\Map(T,X))^T.\label{eq:sqfam}
\eeq
A choice of trivialization of~$\Hodge^{2/2}$ identifies the above with a 2-periodic complex, which is the local value of the sheaf on~$\Mell$ defined in~\S\ref{sec:double}. Under this correspondence, the Bott class $\beta$ is identified with the trivializing section of $\Hodge^{2/2}$. Such local trivializations fail to assemble to a global trivialization of~$\Hodge$; their transformation properties under isomorphisms of elliptic curves leads to a nontrivial sheaf on~$\Mell$. 

\section{Complex analytic elliptic curves and elliptic cohomology}\label{sec:complexanal}

An \emph{elliptic cohomology theory} consists of the data of (1) an elliptic curve $E$ over a commutative ring $R$, (2) an even periodic multiplicative cohomology theory $h$, and (3) isomorphisms $h^0(\pt)\simeq R$ of rings and isomorphisms ${\rm Spf}(h^0(\CP^\infty))\simeq \widehat{E}$  of formal groups over~$R$. Here $\widehat{E}$ is the completion of the elliptic curve~$E$ along its identity section, and ${\rm Spf}(h^0(\CP^\infty))$ is a formal group using naturality of $h$ for the map $\CP^\infty\times\CP^\infty\to \CP^\infty$ classifying the universal tensor product of line bundles. For a comprehensive introduction, see~\cite[\S1]{Lurie_Elliptic}. 

Over the complex numbers, an elliptic curve takes the form $\C/\Lambda$ for $\Lambda\subset \C$ a lattice. The map $\C\to \C/\Lambda$ gives an isomorphism $t_\Lambda$ between the additive formal group law and the group law of~$\C/\Lambda$. For each lattice $\Lambda$ we obtain an elliptic cohomology theory given by the data $(\C/\Lambda,\H(-;\C[\beta,\beta^{-1}]),t_\Lambda)$, i.e., (1) the curve $\C/\Lambda$, (2) ordinary cohomology with values in the graded ring~$\C[\beta,\beta^{-1}]$ with $|\beta|=-2$, and (3) the isomorphism~$t_\Lambda$ (see~\cite[Example~2.3]{AHS}). These data are indeed compatible because the Chern class of a tensor product of line bundles in ordinary cohomology is simply the addition of their Chern classes, i.e., the additive group law from~$(\C,+)$. 

In this section we give a variant of the above for complex analytic families of elliptic curves, and we construct a cocycle refinement to build a sheaf of cdgas over the moduli stack of elliptic curves. Over~$\C$, this plays the role of the sheaf of $E_\infty$-ring spectra whose global sections are~$\TMF$. We review the basic constructions in~\S\ref{sec:defnsheaf1}, compute the derived global sections of the sheaf over the moduli stack of complex analytic curves in~\S\ref{sec:tmfc}, and give a cocycle representative of the derived global section defined by the Witten class of a smooth manifold with a choice of rational string structure in~\S\ref{sec:Wittendif}. 




\subsection{Complex analytic elliptic curves}
The upper half plane $\cH\subset \C$ parameterizes elliptic curves with a specified basis for their 1st homology via the $\Z^2$ quotient 
\beq
&&E=(\cH\times \C)/\Z^2\qquad (\tau,z)\mapsto (\tau,z+m+n\tau),\ \tau\in \cH, z\in \C, (m,n)\in \Z^2,\label{eq:refE}
\eeq
so that the fiber $E_\tau$ at $\tau\in \cH$ is the torus $\C/(\Z\oplus\tau\Z)$. The family~$E$ has an $\SL_2(\Z)$-action 
$$
\tau\mapsto \frac{a\tau+b}{c\tau+d},\quad z\mapsto \frac{z}{c\tau+d},\qquad \left[\begin{array}{cc} a & b \\ c & d\end{array}\right]\in \SL_2(\Z).
$$
Define the \emph{moduli of elliptic curves} $\Mell=\cH/\SL_2(\Z)$ and the \emph{universal elliptic curve} $\EE=E/\SL_2(\Z)$ as holomorphic orbifolds (equivalently, etal\'e groupoids in complex manifolds). The \emph{Hodge bundle} $\Hodge$ is the line bundle on $\Mell$ whose fiber at $E_\tau$ is the space of holomorphic 1-forms on~$E_\tau$. There is a global section of $\Hodge$ over~$\cH$, namely the holomorphic 1-form $\omega$ on $E$ descending from $dz$ on~$\C$. This gives a description of global sections of~$\Hodge^{\otimes k}$
$$
\Gamma(\Mell;\Hodge^{\otimes k})=\left\{f\in \mathcal{O}(\cH)\mid f\left(\frac{a\tau+b}{c\tau+d}\right)=(c\tau+d)^kf(\tau)\right\},
$$
which is the space of \emph{weight $k$ weak modular forms}. Global sections that are also meromorphic as $\tau\to i\infty$ are \emph{weight $k$ weakly holomorphic modular forms,} denoted $\MF_k$. We regard $\MF$ as a graded commutative ring with $\MF^{2k}=\MF_{-k}$ and $\MF^{2k+1}=0$. 

\begin{convention}\label{not:mero}
Hereafter, we impose a regularity condition at infinity for global sections of holomorphic vector bundles over $\Mell$. Specifically, we require sections have meromorphic dependence as~$\tau\to i\infty$ for $\tau\in \cH$. For $\mathcal{V}$ a holomorphic vector bundle, the notation $\Gamma(\Mell,\mathcal{V})$ will always mean global sections with this condition imposed. In particular, $\Gamma(\Mell;\Hodge^{\otimes k})$ denotes weight $k$ weakly holomorphic modular forms. 
\end{convention}

\subsection{Complex analytic elliptic cohomology}\label{sec:defnsheaf1} 

In view of the definition of elliptic cohomology sketched at the beginning of the section, we adopt the following. 

\begin{defn} A \emph{complex analytic elliptic cohomology theory} is (1) a complex analytic family of elliptic curves $E$ over a complex manifold $B$, (2) an even periodic multiplicative cohomology theory $h$, and (3) isomorphisms $h^0(\pt)\simeq \mathcal{O}(B)$ of $\C$-algebras and isomorphisms ${\rm Spf}(h^0(\CP^\infty))\simeq \widehat{E}$ of formal groups over~$\mathcal{O}(B)$. 
\end{defn}


Our primary example of a complex analytic elliptic cohomology theory is the triple given by (1) the complex analytic family of elliptic curves~$E\to \cH$ from~\eqref{eq:refE}, (2) ordinary cohomology valued in the graded ring $\mathcal{O}(\cH)[\beta,\beta^{-1}]$ with $|\beta|=-2$, and (3) the $\cH$-family of holomorphic homomorphisms~$\cH\times \C\to E$ that on fibers induces an isomorphism between the formal additive group and the formal group of $E_\tau$. 
One can put an $\SL_2(\Z)$-action on this cohomology theory to build an equivariant sheaf of cohomology theories over~$\Mell\simeq\cH/\SL_2(\Z)$. However, both the applications to supersymmetric localization and the general philosophy regarding sheaves of cohomology theories on moduli spaces (i.e., the Hopkins--Miller Theorem~\cite[Theorem~1.1]{Lurie_Elliptic}) require that we refine this object. Namely, we pass from a sheaf of cohomology groups to a sheaf of commutative differential graded algebras on~$\cH/\SL_2(\Z)\simeq \Mell$. We implement this as follows.

\begin{defn} \label{defn:Ell1}
For a smooth manifold $X$, define a sheaf $(\Ell(X),d)$ on $\cH$ of commutative differential graded algebras that to $U\subset\cH$ assigns the complex of $\mathcal{O}(U)$-modules
\beq
&&(\Ell(X)(U),d):= \Big(\mathcal{O}(U;\Omega^\bullet (X)[\beta,\beta^{-1}]),d\Big)\simeq \Big(\Omega^\bullet (X;\mathcal{O}(U)[\beta,\beta^{-1}]),d\Big) \nonumber 
\eeq
where $d$ is the de~Rham differential and the grading is the total grading for a differential form valued in the graded ring $\mathcal{O}(U)[\beta,\beta^{-1}]$ with $|\beta|=-2$. We promote this to an $\SL_2(\Z)$-equivariant sheaf that to a map $U\to U'\subset \cH$ coming from $A \in \SL_2(\Z)$ assigns the map
$$
f_A\colon \Ell(X)(U')\to \Ell(X)(U)\qquad A=\left[\begin{array}{cc} a & b \\ c & d\end{array}\right]\in \SL_2(\Z)
$$
that first pulls back $\Omega^\bullet (X)[\beta,\beta^{-1}]$-valued functions, and then applies the algebra automorphism determined by~$\beta\mapsto \beta /(c\tau+d)$. This $\SL_2(\Z)$-equivariant sheaf on $\cH$ determines a sheaf on $\Mell=\cH/\SL_2(\Z)$. 
%
\end{defn}

\begin{rmk}
With its focus on smooth manifolds and differential forms, the above can be thought of as a version of \emph{differential} complex analytic elliptic cohomology. One could instead take singular cochains on topological spaces with values in~$\mathcal{O}(U)[\beta,\beta^{-1}]$. The focus on manifolds is for our intended contact with ideas from localization in physics, e.g., equivariant de~Rham complexes of double loop spaces. 
\end{rmk}

\begin{rmk}
One can view $(\Ell(X)(U),d)$ as the de~Rham cohomology of $X$ with values in powers of the Hodge bundle, $\Hodge^{\otimes \bullet}$. The class $\beta^k$ is a trivialization of $\Hodge^{\otimes k}$ on $\HH$, identifying sections with functions $\mathcal{O}(\HH)\cdot \beta^k\simeq \mathcal{O}(\HH;\Hodge^{\otimes k})$.
\end{rmk}

\subsection{Complex analytic topological modular forms} \label{sec:tmfc}

For a suborbifold $U/\Gamma\subset \cH/\SL_2(\Z)$, taking sections gives a functor from $\Gamma$-equivariant complexes of sheaves to chain complexes. Applied to the restriction of $\Ell(X)$ to $U/\Gamma$, we have
\beq
(\Ell(X),d)|_{U/\Gamma}\mapsto (\Ell(X)(U)^\Gamma,d)=\left(\mathcal{O}\Big(U;\Omega^\bullet (X)[\beta,\beta^{-1}]\right)^\Gamma,d\Big).\label{eq:functor}
\eeq

\begin{defn} The right derived functors of~\eqref{eq:functor} are the \emph{hypercohomology} of~$(\Ell(X),d)$ on~$U/\Gamma$, denoted by $\HH^\bullet(U/\Gamma;\Ell(X))$. \end{defn}


\begin{proof}[Proof of Proposition~\ref{prop:1}]
Part 1 of the proposition follows from the definition. In proving part~2, we will set up the general description of the hypercohomology for an arbitrary suborbifold~$U/\Gamma$, specializing to the case in the statement at the end of the proof (this setup is for later convenience). To compute the hypercohomology, we use the Dolbeault resolution of holomorphic sections with differential $\overline{\partial}$, and the standard resolution from group cohomology of a $\Gamma$-module using the differential~$\delta$ that takes alternating sums of pullbacks along the action by~$\Gamma^{\times k}$. Along with the differential~$Q$ on $\Ell(X)$, this gives a triple complex computing the hypercohomology. For a fixed degree~$k$ in $\Ell(X)$, we have the (doubly augmented) double complex
\beq
\begin{tikzpicture}[baseline=(basepoint)];
\node (A0) at (-10,0) {$\mathcal{O}(U;\mathcal{E}^k)^{\Gamma}$}; 
\node (A00) at (-8.5,-.5) {$\null$};
\node (A02) at (-10.5,.5) {$\null$};
\node (A1) at (-6.75,0) {$\Omega^{0,0}(U;\mathcal{E}^k)^\Gamma$}; 
\node (A2) at (-3,0) {$\Omega^{0,1}(U;\mathcal{E}^k)^\Gamma$}; 
\node (A3) at (0,0) {$0$};
\node (B0) at (-10,1) {$\mathcal{O}(U;\mathcal{E}^k)$}; 
\node (B1) at (-6.75,1) {$\Omega^{0,0}(U;\mathcal{E}^k)$}; 
\node (B2) at (-3,1) {$\Omega^{0,1}(U;\mathcal{E}^k)$}; 
\node (B3) at (0,1) {$0$};
\node (C0) at (-10,2) {$\mathcal{O}(U\times \Gamma;\mathcal{E}^k)$}; 
\node (C1) at (-6.75,2) {$\Omega^{0,0}(U\times \Gamma;\mathcal{E}^k)$}; 
\node (C2) at (-3,2) {$\Omega^{0,1}(U\times \Gamma;\mathcal{E}^k)$}; 
\node (C3) at (0,2) {$0$};
\node (D0) at (-10,3) {$\mathcal{O}(U\times \Gamma^{\times 2};\mathcal{E}^k)$}; 
\node (D1) at (-6.75,3) {$\Omega^{0,0}(U\times \Gamma^{\times 2};\mathcal{E}^k)$}; 
\node (D2) at (-3,3) {$\Omega^{0,1}(U\times \Gamma^{\times 2};\mathcal{E}^k)$}; 
\node (D3) at (0,3) {$0$};
\node (E0) at (-10,4) {$\vdots$}; 
\node (E00) at (-8.5,4) {$\null$};
\node (E1) at (-6.75,4) {$\vdots$}; 
\node (E2) at (-3,4) {$\vdots$}; 
\node (E22) at (-1,.5) {$\null$};
\draw[->] (A0) to (A1);
\draw[->] (A1) to node [above] {$\overline{\partial}$} (A2);
\draw[->] (A2) to (A3);
\draw[->] (B0) to (B1);
\draw[->] (B1) to node [above] {$\overline{\partial}$} (B2);
\draw[->] (B2) to (B3);
\draw[->] (C0) to (C1);
\draw[->] (C1) to node [above] {$\overline{\partial}$} (C2);
\draw[->] (C2) to (C3);
\draw[->] (D0) to (D1);
\draw[->] (D1) to node [above] {$\overline{\partial}$} (D2);
\draw[->] (D2) to (D3);
\draw[->] (A0) to (B0);
\draw[->] (B0) to node [left] {$\delta$} (C0);
\draw[->] (C0) to node [left] {$\delta$} (D0);
\draw[->] (A1) to (B1);
\draw[->] (B1) to node [left] {$\delta$} (C1);
\draw[->] (C1) to node [left] {$\delta$} (D1);
\draw[->] (A2) to (B2);
\draw[->] (B2) to node [left] {$\delta$} (C2);
\draw[->] (C2) to node [left] {$\delta$} (D2);
\draw[->] (D0) to node [left] {$\delta$} (E0);
\draw[->] (D1) to node [left] {$\delta$} (E1);
\draw[->] (D2) to node [left] {$\delta$} (E2);
\draw[-,dashed] (A00) to (E00);
\draw[-,dashed] (A02) to (E22);
\path (0,1.5) coordinate (basepoint);
\end{tikzpicture}\label{diagram:double}
\eeq
where $\mathcal{E}^k$ denotes the subsheaf of degree~$k$ sections of~$\Ell(X)$.
We compute hypercohomology by first taking the cohomology of these double complexes for each $k\in \Z$. 

We claim the rows of~\eqref{diagram:double} are exact. Indeed, the $\overline{\partial}$-cohomology is concentrated in degree zero because any open subset~$U\subset \cH$ has Dolbeault cohomology concentrated in degree zero (e.g., because it is a Stein manifold); sections of $\Ell(X)$ over~$U$ can be identified with sections of the trivial holomorphic vector bundle with fiber $\Omega^\bullet(X)[\beta,\beta^{-1}]$, so Dolbeault cohomology valued in this trivial bundle is also concentrated in degree zero. The columns of~\eqref{diagram:double} are also exact: the cohomology with respect to~$\delta$ is the groupoid cohomology of $U/\Gamma$ with values in the $\Gamma$-equivariant vector bundle on $U$ of $(0,j)$-forms valued in the trivial bundle with fiber $\Omega^\bullet(X)[\beta,\beta^{-1}]$. This is known to be concentrated in degree zero for proper groupoids; for example, see~\cite[Proposition~1]{Crainic}. We now apply a standard argument from homological algebra: if the rows of an augmented double complex are exact, then the total cohomology is isomorphic to the cohomology of the initial column.
 This implies that the cohomology of the first column (excluding the first term) equals the cohomology of the first row (excluding the first term), and both of these are equal to the cohomology of the double complex. 

Specializing to part (2) of the proposition, we use that the orbifold $\cH/\SL_2(\Z)$ is Morita equivalent to $(\C\setminus \{0,1\})/(\Z/2\times S_3)$ where $S_3$ is the symmetric group on~3 letters. Equivalently, the complex analytic stacks underlying these orbifolds are isomorphic. See, for example,~\cite[Proposition~3.15]{hain} for a construction of this isomorphism. Groupoid cohomology is invariant under Morita equivalences~\cite[\S3]{BehrendXu}. This allows us to replace the first column of the double complex~\eqref{diagram:double} with the analogous one for $(\C\setminus \{0,1\})/(\Z/2\times S_3)$ which computes the group cohomology of $\Z/2\times S_3$ with coefficients in the module $\mathcal{O}(\C\setminus \{0,1\};\mathcal{E}^k)$. But the order of this finite group is invertible in the module, so the cohomology above degree zero vanishes. Hence the cohomology of the double complex~\eqref{diagram:double} is~$\mathcal{O}(\cH;\mathcal{E}^k)^{\SL_2(\Z)}$, i.e., the $\SL_2(\Z)$-invariant holomorphic sections. This identifies the hypercohomology with cohomology of the de~Rham complex of~$X$ with values in the graded ring~$\MF$, i.e.,~$\H(X;\MF)\simeq \TMF(X)\otimes \C$. This completes the proof. \ep

\subsection{A cocycle representative of the Witten class as a derived global section}\label{sec:Wittendif}
The \emph{Witten class} is the characteristic class
\beq
&&[\Wit(X)]:=\exp\left( \sum_{k=1}^\infty [{ph}_{4k}(TX)] \beta^{2k}E_{2k}\right)\in \H^\bullet(X;\mathcal{O}(\cH)[\beta,\beta^{-1}]).\label{eq:WC}
\eeq
where $[{ph}_{4k}(TX)]\in \H^{4k}(TX)$ is the $4k$th component of the Pontryagin character and $E_{2k}\in \mathcal{O}(\cH)$ is the $2k$th Eistenstein series, 
$$
E_{2k}(\tau)=\sum_{n,m\in \Z_*^2} \frac{1}{(m\tau+n)^{2k}},\quad \qquad \Z_*^2=\{(n,m)\in \Z^2 \mid (n,m)\ne (0,0)\}
$$ 
and we take~$E_2$ to be the holomorphic version of the 2nd Eisenstein series (the above sum is only conditionally convergent when $k=1$). Since $E_{2k}$ is a modular form for $k>1$, we have $[\Wit(X)]\in \H^\bullet(X;\MF)\simeq \TMF(X)\otimes \C$ if and only if $[{ph}_1(TX)]=[p_1(TX)]=0$. In this case, the \emph{Witten genus} 
$$
\int_X [\Wit(X)]\in \MF^{-{\rm dim}(X)}=\MF_{{\rm dim}(X)/2}
$$
is a modular form of weight ${\rm dim}(X)/2$. 

\begin{proof}[Proof of Proposition~\ref{prop:2}]
For a choice of Riemannian metric on $X$, let $R\in \Omega^2(X;\End(TM))$ denote the curvature of the Levi-Civita connection and $ph_{2k}(TM)={\rm Tr}(R^{2k})/(2k(2\pi i )^{2k})$ be a cocycle representative for the degree~$4k$ part of the Pontryagin character. Consider the cocycle refinement of the Witten class,
$$
\Wit(X):=\exp\left( \sum_{k=1}^\infty \beta^{2k} {ph}_{2k}(TM)E_{2k}\right)\in \mathcal{O}(\cH;\Omega^\bullet(X)[\beta,\beta^{-1}]),
$$
but note that it fails to descend to a derived global section because $E_2$ is not a modular form. Instead,
\beq
E_2\left(\frac{a\tau+b}{c\tau+d}\right)=(c\tau+d)^2E_2(\tau)-2\pi i c(c\tau+d),\qquad \left[\begin{array}{cc} a & b \\ c & d\end{array}\right]\in \SL_2(\Z).\label{eq:E2}
\eeq
If $[p_1(TX)]=0$, then $\delta\Wit(X)$ is $d$-exact (with $\delta$ as in~\eqref{diagram:double}). A rational string structure on~$X$---i.e., $H\in \Omega^3(X)$ with $dH=p_1(TX)$---uniquely determines an element $A \in \mathcal{O}(\cH\times \SL_2(\Z);\Omega^\bullet(X)[\beta,\beta^{-1}])$ satisfying
\beq
&&\delta\Wit(X)=dA \qquad A:=H\left(\Wit(X)\frac{1- \exp(-p_1(TX) 2\pi i c/(c\tau+d))}{p_1(TX)}\right) \label{eq:Witdelta}
\eeq
where we use that $1- \exp(-p_1 2\pi i c/(c\tau+d))$ is divisible by $p_1(TX)$ in the formula for~$A$. 
Therefore $(\Wit(X),A)$ defines a degree zero cocycle in~\eqref{diagram:double}, and hence a derived global section whose underlying class in $\TMF^0(X)\otimes \C$ is the Witten class of $X$. 
\ep


%

\section{Localization and a double loop space model for elliptic cohomology} \label{sec:double}

\subsection{More geometry of elliptic curves over $\C$}\label{sec:EC}

Endow the family $E\to \cH$ from~\eqref{eq:refE} with the fiberwise flat metric descending from~$\C$, normalized so that the fibers have unit volume. Recall that~$\omega$ is the 1-form on $E$ descending from~$dz$ on~$\cH\times \C$. Consider the following vector field on~$E$. 

\begin{defn}
Let~$\omega^\sharp$ be the vector field dual to~$\omega$ via the K\"ahler form on~$E$. 
\end{defn}

Since the K\"ahler form is $\SL_2(\Z)$-invariant, $\omega^\sharp$ transforms under $\SL_2(\Z)$ as
\beq
&&\omega^\sharp\mapsto \frac{1}{c\tau+d}\omega^\sharp.\label{eq:omegatransf}
\eeq
For future use, we note the diffeomorphism from $E$ to the constant family of square tori
\beq
&&\cH\times (\R^2/\Z^2)\stackrel{\sim}{\to} (\cH\times \C)/\Z^2=E, \qquad (\tau,x,y)\mapsto (\tau,x+\tau y)\quad (x,y)\in \R^2.\label{eq:square}
\eeq
This is $\SL_2(\Z)$-equivariant for the $\SL_2(\Z)$-action on $\R^2$ given by~$x\mapsto ax-by, y\mapsto dy-cx.$ The pullback of the holomorphic 1-form $\omega$ along~\eqref{eq:square} is~$dx+\tau dy,$ and the K\"ahler form is $dx\wedge dy$. Hence $\omega^\sharp$ in this description is $\tau\partial_x-\partial_y$.
We also have $\overline{\omega}^\sharp=\overline{\tau}\partial_x-\partial_y$ for the vector field $\overline{\omega}^\sharp$ dual to $\overline{\omega}$.

\subsection{The sheaf $(\Ell(\X),Q)$ on $\Mell$}\label{sec:defnsheaf}
Let $\X\in \Sh(\TMfld)$ be a sheaf on the site of $T$-manifolds, and define $T$-invariant differential forms on $\X$, denoted $\Omega^\bullet(\X)^T$, as morphisms of sheaves $\X\to \Omega^\bullet(-)^T$; see~\S\ref{sec:appen} for details. The main example will be $T$-invariant forms on~$\X=\Map(T,X)$. The isomorphism~\eqref{eq:square} determines a complex vector field~$\omega^\sharp=\tau\partial_x-\partial_y$ on~$\cH\times T$. Let $\iota_{\omega^\sharp}=\tau\iota_{\partial_x}-\iota_{\partial_y}$ be the associated contraction operator on $\mathcal{O}(\cH;\Omega^\bullet(\X))\simeq \Omega^\bullet(\X;\mathcal{O}(\HH))$ for~$\X\in \Sh(\TMfld)$. 

\begin{defn} \label{defn:Ell}
For $\X\in \Sh(\TMfld)$, define a sheaf $(\Ell(\X),Q)$ on $\cH$ of commutative differential graded algebras that to $U\subset\cH$ assigns the complex of $\mathcal{O}(U)$-modules
\beq
&&(\Ell(\X)(U),Q):= \left(\mathcal{O}\Big(U;\Omega^\bullet (\X;\C[\beta,\beta^{-1}])^T\right),Q\Big)\label{eq:Ell}
\eeq
where $Q=d-\beta^{-1}\iota_{\omega^\sharp}$ and $|\beta|=-2$. We promote this to an $\SL_2(\Z)$-equivariant sheaf that to a map $U\to U'\subset \cH$ coming from $A \in \SL_2(\Z)$, assigns the map
$$
f_A\colon \Ell(\X)(U')\to \Ell(\X)(U)\qquad A=\left[\begin{array}{cc} a & b \\ c & d\end{array}\right]\in \SL_2(\Z)
$$
that first pulls back $\Omega^\bullet (\X;\C[\beta,\beta^{-1}])^T$-valued functions, and then applies the algebra automorphism determined by~$\beta\mapsto \beta /(c\tau+d)$. This $\SL_2(\Z)$-equivariant sheaf on $\cH$ determines a sheaf on $\Mell=\cH/\SL_2(\Z)$. 
%
\end{defn}

\begin{rmk} When $\X=X$ is a representable presheaf with trivial $T$-action, the sheaf $(\Ell(X),d)$ from Definition~\ref{defn:Ell1} agrees with Definition~\ref{defn:Ell}: the trivial action implies $Q=d$.
\end{rmk}

\begin{rmk} As usual, the element $\beta$ in~\eqref{eq:Ell} can be identified with the choice of invariant 1-form~$\omega$ on the curve. This gives a geometric reason for why the differential~$Q$ is compatible with the $\SL_2(\Z)$-equivariant structure. Algebraically this compatibility is because~$\omega^\sharp$ and $\beta$ transform the same way under $\SL_2(\Z)$; see~\eqref{eq:omegatransf}. 
\end{rmk}

\begin{rmk} \label{rmk:power}
In the sheaf-theoretic definition of differential forms on $\X$ adopted in this paper, the algebra~$\Omega^\bullet (\X;\C[\beta,\beta^{-1}])^T$ implements a completion in powers of~$\beta$. Recall that a map $\omega\colon \X\to \Omega^k (-;\C[\beta,\beta^{-1}])^T$ assigns to each $S$-point of $\X$ a differential form $\omega_S\in \Omega^k (S;\C[\beta,\beta^{-1}])^T$, and these are compatible under base change. As we vary~$S$, the polynomial dependence in~$\beta$ can have unbounded degree, e.g., on a sequence of $S_i$-points of $\X$ where $S_i$ has unbounded dimension. This feature is crucial in localization arguments; see Remark~\ref{rmk:powerbeta}. More formally,~$\X$ is a colimit of its $S$-points, so $\Omega^\bullet (\X;\C[\beta,\beta^{-1}])^T$ is a limit of graded rings $\Omega^\bullet (S;\C[\beta,\beta^{-1}])^T$. For this reason, one might loosely regard~$\Omega^\bullet (-;\C[\beta,\beta^{-1}])^T$ as~$\Omega^k (\X;\C(\!(\beta)\!))^T$. We avoid the latter notation simply because $\C(\!(\beta)\!)$ is not a graded ring with~$|\beta|=-2$.
\end{rmk}

\begin{rmk} We observe that $\Ell(\X)(U)$ is the sheaf of holomorphic sections of a trivial bundle on~$\cH$ with fiber the graded vector space $\Omega^\bullet (\X;\C[\beta,\beta^{-1}])^T$. The $\SL_2(\Z)$-equivariant structure for the sheaf $\Ell(\X)(U)$ is inherited through the $\SL_2(\Z)$-action on~$\beta$. Furthermore,~$Q$ comes from a degree~$+1$ vector bundle automorphism. \end{rmk}

From the final remark above, the following is a straightforward consequence of naturality of~$\Omega^\bullet (\X;\C[\beta,\beta^{-1}])^T$ in~$\X$ and properties the sheaf $\Omega(-)^T$ on $T$-manifolds. 

\begin{lem} \label{lem:homotopy}A map $f\colon \X\to \Y$ in $\Sh(\TMfld)$ induces a morphism of sheaves of cdgas $f^*\colon (\Ell(\Y),Q)\to (\Ell(\X),Q)$. Smoothly homotopic maps
 induce chain homotopic maps of sheaves of cdgas. \end{lem}

\subsection{Localization}

There is a monomorphism of sheaves $X\hookrightarrow \Map(T,X)$ including~$X$ as the constant maps, or equivalently, the $T$-fixed subsheaf. Define the sheaf $\Map(T,X)\setminus X$ as having $S$-points maps $S\to \Map(T,X)$ such that the pullback $X\hookrightarrow \Map(T,X)\leftarrow S$ is empty. 

\begin{lem} 
There exists an element $\alpha\in C^\infty(\cH;\Omega^1(\Map(T,X)\setminus X)^T)$ such that~$\iota_{\omega^\sharp}\alpha=1.$
\end{lem}

\begin{proof} Choose a $T$-invariant Riemannian metric on $X$. For a $T$-invariant map $\phi\colon S\times T\to X$ associated to an $S$-point of $\Map(T,X)\setminus X$, consider the 1-form on $S$
\beq
\alpha_S:= \frac{\int_T\langle \overline{\omega}^\sharp,-\rangle}{\int_T\langle \overline{\omega}^\sharp,\omega^\sharp\rangle}\in C^\infty(\cH;\Omega^1(S))\simeq \Omega^1(S;C^\infty(\cH))\label{eq:alphadef}
\eeq
where $\langle-,-\rangle$ denotes the $\C$-bilinear extension of the pullback of the Riemannian metric along $\phi$, and $\omega^\sharp,\overline{\omega}^\sharp$ denote the complex vector fields on $\cH\times S\times T$ coming from the infinitesimal $T$-action. To verify that $\alpha$ is well-defined, we first calculate (using formulas from the end of~\S\ref{sec:EC})
\beq
\langle \overline{\omega}^\sharp,\omega^\sharp\rangle&=&\langle \overline{\tau}\partial_x-\partial_y,\tau\partial_x-\partial_y\rangle=\langle (\tau_0\partial_x-\partial_y)-i\tau_1\partial_x,(\tau_0\partial_x-\partial_y)+i\tau_1\partial_x\rangle\nonumber\\
&=&|\tau_0\partial_x-\partial_y|^2+|\tau_1\partial_x|^2\in C^\infty(\cH\times S\times T)\nonumber 
\eeq
where $\tau=\tau_0+i\tau_1\in \cH$ for $\tau_0,\tau_1\in \R$ with $\tau_1>0$. Next observe that for a family of nonconstant maps, if $\int_T |\partial_x|^2=0$, then $\int_T|\partial_y|^2> 0$. Therefore, for an $S$-point of $\Map(T,X)\setminus X$ we have that $\int_T\langle \overline{\omega}^\sharp,\omega^\sharp\rangle\in C^\infty(\cH\times S;\R_{>0})$ is real and positive. Furthermore, $\alpha_S$ is $T$-invariant by construction and satisfies $\iota_{\omega^\sharp}\alpha_S=1$. Finally, $\alpha_S$ is natural in $S$, so we obtain a 1-form $\alpha\in C^\infty(\cH;\Omega^1(\Map(T,X)\setminus X)^T)$ with the desired property. 
\ep

From~\eqref{eq:omegatransf} and~\eqref{eq:alphadef}, we observe that $\beta\alpha$ is $\SL_2(\Z)$-invariant. Furthermore, if $\eta\in \Omega^{0,j}(\cH;\Omega^\bullet(\Map(T,X)\setminus X;\C[\beta,\beta^{-1}])^T)^{\SL_2(\Z)}$ is $Q$-closed, then it is $Q$-exact via 
\beq
Q\left(\frac{\beta \alpha}{\beta d\alpha-1}\eta\right)=Q\left(\frac{\beta\alpha}{\beta d\alpha-1}\right)\eta=\frac{\beta d\alpha-\iota_{\omega^\sharp}\alpha}{\beta d\alpha-1}\eta=\eta.\label{eq:powerbeta}
\eeq
This leads one to expect the cohomology of $\Ell(\Map(T,X))$ to be supported in a neighborhood of $X\subset \Map(T,X)$, as the following lemma verifies.

\begin{lem}\label{lem:openrest} Let $N\subset \Map(T,X)$ be an open subsheaf and $\part\in C^\infty(\Map(T,X))^T$ such that
\begin{enumerate}
\item  The constant maps are a subsheaf of $N$, i.e., there is a factorization $X\subset N\subset \Map(T,X)$ of the inclusion $X\subset \Map(T,X)$;
\item $\part\equiv 1$ on an open subsheaf containing~$X$ and the support of $\part$ is contained in $N$. 
\end{enumerate} 
Then the restriction map $(\Ell(\Map(T,X)),Q) \to (\Ell(N),Q)$ induces an isomorphism on hypercohomology for each suborbifold $U/\Gamma\subset \cH/\SL_2(\Z)$.
\end{lem}
\bp 
Consider the double complex~\eqref{diagram:double} where now $\mathcal{E}^k$ is the degree $k$ subsheaf of $\Ell(\Map(T,X))$, rather than of~$\Ell(X)$. As before, the rows and columns of this double complex are exact, so the cohomology can be computed using the first row. Hence, to demonstrate a quasi-isomorphism on hypercohomology, it suffices to show that the restriction map for $N\hookrightarrow \Map(T,X)$ induces an isomorphism on $Q$-cohomology in the columns to the right of the dashed line:
\beq
\begin{tikzpicture}[baseline=(basepoint)];
\node (A0) at (-10,0) {$\mathcal{O}(U;\mathcal{E}^0)^{\Gamma}$}; 
\node (A1) at (-6.75,0) {$\Omega^{0,0}(U;\mathcal{E}^0)^\Gamma$}; 
\node (A2) at (-3,0) {$\Omega^{0,1}(U;\mathcal{E}^0)^\Gamma$}; 
\node (A3) at (0,0) {$0$};
\node (A00) at (-8.5,-0.5) {$\null$};
\node (B0) at (-10,1) {$\mathcal{O}(U;\mathcal{E}^1)^\Gamma$}; 
\node (B1) at (-6.75,1) {$\Omega^{0,0}(U;\mathcal{E}^1)^\Gamma$}; 
\node (B2) at (-3,1) {$\Omega^{0,1}(U;\mathcal{E}^1)^\Gamma$}; 
\node (B3) at (0,1) {$0$};
\node (C0) at (-10,2) {$\mathcal{O}(U;\mathcal{E}^2)^\Gamma$}; 
\node (C1) at (-6.75,2) {$\Omega^{0,0}(U;\mathcal{E}^2)^\Gamma$}; 
\node (C2) at (-3,2) {$\Omega^{0,1}(U;\mathcal{E}^2)^\Gamma$}; 
\node (C3) at (0,2) {$0$};
\node (D0) at (-10,3) {$\mathcal{O}(U;\mathcal{E}^3)^\Gamma$}; 
\node (D1) at (-6.75,3) {$\Omega^{0,0}(U;\mathcal{E}^3)^\Gamma$}; 
\node (D2) at (-3,3) {$\Omega^{0,1}(U;\mathcal{E}^3)^\Gamma$}; 
\node (D3) at (0,3) {$0$};
\node (E0) at (-10,4) {$\vdots$}; 
\node (E1) at (-6.75,4) {$\vdots$}; 
\node (E2) at (-3,4) {$\vdots$}; 
\node (E00) at (-8.5,4) {$\null$};
\draw[->] (A0) to (A1);
\draw[->] (A1) to node [above] {$\overline{\partial}$} (A2);
\draw[->] (A2) to (A3);
\draw[->] (B0) to (B1);
\draw[->] (B1) to node [above] {$\overline{\partial}$} (B2);
\draw[->] (B2) to (B3);
\draw[->] (C0) to (C1);
\draw[->] (C1) to node [above] {$\overline{\partial}$} (C2);
\draw[->] (C2) to (C3);
\draw[->] (D0) to (D1);
\draw[->] (D1) to node [above] {$\overline{\partial}$} (D2);
\draw[->] (D2) to (D3);
\draw[->] (A0) to node [left] {$Q$} (B0);
\draw[->] (B0) to node [left] {$Q$} (C0);
\draw[->] (C0) to node [left] {$Q$} (D0);
\draw[->] (A1) to node [left] {$Q$} (B1);
\draw[->] (B1) to node [left] {$Q$} (C1);
\draw[->] (C1) to node [left] {$Q$} (D1);
\draw[->] (A2) to node [left] {$Q$} (B2);
\draw[->] (B2) to node [left] {$Q$} (C2);
\draw[->] (C2) to node [left] {$Q$} (D2);
\draw[->] (D0) to node [left] {$Q$} (E0);
\draw[->] (D1) to node [left] {$Q$} (E1);
\draw[->] (D2) to node [left] {$Q$} (E2);
\draw[-,dashed] (A00) to (E00);
\path (0,1.5) coordinate (basepoint);
\end{tikzpicture}\label{diagram:ddouble}
\eeq
These quasi-isomorphisms in $Q$-cohomology mimic the proof of localization in finite dimensions \cite[Proposition~2.1]{BerlineVergne2}. 
First we show that the restriction map is surjective on $Q$-cohomology. So given $\mu\in \Omega^{0,j}(U;\Omega^\bullet(N;\C[\beta,\beta^{-1}])^T)^\Gamma$ with $Q\mu=0$, consider
 $\widetilde{\mu}\in \Omega^{0,j}(U;\Omega^\bullet(\Map(T,X);\C[\beta,\beta^{-1}])^T)^\Gamma$ the $Q$-closed section defined 
 by
\beq
\widetilde{\mu}=\left\{\begin{array}{ll} \mu-Q\left((1-\part)\frac{\beta \alpha}{\beta d\alpha-1}\mu\right) & {\rm on} \ N\setminus X \subset \Map(T,X)\setminus X\\ \mu & {\rm on}\ X\\ 
0 & {\rm otherwise}, \end{array}\right.\label{eq:mutil1}
\eeq
where $\alpha\in C^\infty(U;\Omega^1(\Map(T,X)\setminus X)^T)$ is from the previous lemma and $\beta\alpha$ is $\Gamma$-invariant. This is well-defined by the observation immediately before the statement of the present lemma, together with the fact that $\part\equiv 1$ on an open subsheaf of $N$ containing~$X$. Finally, $\widetilde{\mu}$ is $Q$-cohomologous to $\mu$ on $N$, so the restriction is surjective on $Q$-cohomology. 

Next we show that the restriction is injective on $Q$-cohomology. Given a section $\mu\in \Omega^{0,j}(U;\Omega^\bullet(\Map(T,X);\C[\beta,\beta^{-1}])^T)^\Gamma$ with $Q\mu=0$ and whose restriction to $N$ satisfies $\mu=Q\gamma$, consider the cohomologous, $Q$-exact section
\beq
\widetilde{\mu}=\left\{ \begin{array}{ll} \mu-Q(\part\gamma)=Q\left(\frac{\beta\alpha}{\beta d\alpha-1}(\mu-Q(\part\gamma))\right) &{\rm on} \ N\setminus X\\ \mu-Q(\part\gamma)=0 & {\rm on}\ X \\
\mu=Q\left(\frac{\beta\alpha}{\beta d\alpha-1}\mu\right) & {\rm otherwise}
\end{array}\right.\label{eq:mutil2}
\eeq
which is $Q$-cohomologous to $\mu$ and is $Q$-exact in $\Omega^{0,j}(U;\Omega^\bullet(\Map(T,X);\C[\beta,\beta^{-1}])^T)^\Gamma$.
Hence, the restriction is also injective on $Q$-cohomology. This shows that restriction gives a quasi-isomorphism of double complexes~\eqref{diagram:ddouble} computing the hypercohomology, and hence the hypercohomology groups are isomorphic. 
\ep

\begin{rmk}\label{rmk:powerbeta} In~\eqref{eq:powerbeta} and similar formulas, we emphasize that the restriction of $(\beta d\alpha-1)^{-1}$ to an $S$-point of $\Map(T,X)$ will only involve finitely many powers of $\beta$. This is because $d\alpha$ has positive degree, so $(d\alpha_S)^N=0$ for large $N$ on any finite-dimensional manifold $S$. 
\end{rmk}

\begin{cor} \label{cor:openrest}
Let $N\subset \Map(T,X)$ be as in the previous lemma. The restriction map $(\Omega^\bullet(\Map(T,X))^{T},Q_\tau)\to (\Omega^\bullet(N)^{T},Q_\tau)$ is a quasi-isomorphism of $\Z/2$-graded complexes. 
\end{cor}
\bp
The cocycles $\widetilde{\mu}\in \Omega^\bullet(\Map(T,X))^T$ defined by the same formulas~\eqref{eq:mutil1} and~\eqref{eq:mutil2} (with $\beta=1$ and $\alpha=\alpha_\tau \in \Omega^1(\Map(T,X)\setminus X)^T$) show that the restriction map on $Q$-cohomology is surjective and injective, respectively. 
\ep

\begin{lem} \label{lem:openexists}There exists an open subsheaf $N\subset \Map(T,X)$ satisfying the conditions in Lemma~\ref{lem:openrest}, that in addition has a deformation retraction to $X\subset N$.\end{lem}

\bp
Fix a metric on $X$ and choose an $\epsilon>0$ smaller than the injectivity radius. Define the sheaf $N$ as having $S$-points
$$
N(S)=\{\phi_0\colon S\to X, v\in \Gamma (S\times T,p^*\phi_0^*TX)\mid \int_T v=0, \|v\|<\epsilon \}
$$
where $p\colon S\times T\to S$ is the projection and the inequality $\|v\|<\epsilon$ is pointwise on $S\times T$, using the pullback of the metric on $X$. Pulling back sections along base changes $S'\to S$ promotes these values for each~$S$ to a sheaf. The Riemannian exponential map applied to~$v$ gives a map~$N\to \Map(T,X)$. This realizes $N$ as the open subsheaf of $\Map(T,X)$ consisting of maps $S\times T\to X$ with the property that for every $s\in S$, the image of $\{s\}\times T\to X$ is contained in an $\epsilon$-ball in~$X$. The factorization $X\subset N\subset\Map(T,X)$ on $S$-points corresponds to the inclusion $X(S)\hookrightarrow N(S)$ along $v=0$. We observe that $N$ has a $T$-equivariant deformation retraction to $X$ where the homotopy~$N\times I\to N$ on $S$-points is given by $(\phi_0,v)\mapsto (\phi_0,(1-t)v)$ for $t\in [0,1]$. 

It remains to construct the function $f$ on $\Map(T,X)$ with desired properties. Consider the function $E$ on $N$ that to an $S$-point $(\phi_0,v)$ assigns
$$
E(\phi_0,v)=\int_T \|v\|^2\in C^\infty(S). 
$$
Choose a bump function $b\colon \R\to \R$ that is 1 in a neighborhood of $0$ and has support in $(-\infty,\epsilon'^2)$ for some positive $\epsilon'<\epsilon$. Then consider the $T$-invariant function on $\Map(T,X)$ whose value on an $S$-point as above is given by
\beq
\part(\phi)=\left\{\begin{array}{cl} b(E) & {\rm on} \ S' \\ 0 & {\rm else}\end{array}\right.\label{eq:bump}
\eeq
where $S'\subset S$ is the open submanifold given by the pullback~$S\to \Map(T,X)\leftarrow N$. To see that this function is well-defined, first observe for any $S$-point of $N$ we have $E(\phi_0,v)<\epsilon$, so $\part$ vanishes on the open subsheaf $E^{-1}((-\epsilon',\epsilon))\subset N\subset \Map(T,X)$. Next we observe that $f(\phi)=0$ on the complement of the closed submanifold $E^{-1}((-\infty,\epsilon'])\subset S'\subset S$. So by the sheaf property, these compatible values on open subsets of $S$ glue to a smooth function on~$S$. Finally, by construction $f$ has support contained in $N$ and $f\equiv 1$ on an open subsheaf containing~$X$. 
\ep

\begin{proof}[Proof of Theorems~\ref{thm:0} and~\ref{thm:1}]
The subsheaf~$N$ from Lemma~\ref{lem:openexists} factors the inclusion of constant maps $X\hookrightarrow N\hookrightarrow \Map(T,X)$, which factors the restriction of sheaves as
$$\Ell(\Map(T,X)) \to \Ell(N)\to \Ell(X),\qquad\qquad \Omega^\bullet(\Map(T,X))^{T}\to \Omega^\bullet(N)^{T}\to  \Omega^\bullet(X)$$ 
in the two respective cases. Lemma~\ref{lem:openrest} shows that $\Ell(\Map(T,X)) \to \Ell(N)$ induces an isomorphism on hypercohomology for any suborbifold of~$\Mell$, and Corollary~\ref{cor:openrest} shows that $(\Omega^\bullet(\Map(T,X))^{T},Q_\tau)\to (\Omega^\bullet(N)^{T},Q_\tau)$ is a quasi-isomorphism. The $T$-equivariant deformation retraction from $N$ to $X$ along with Lemma~\ref{lem:homotopy} shows that the arrow $\Ell(N)\to \Ell(X)$ induces an isomorphism in hypercohomology and that $(\Omega^\bullet(N)^{T},Q_\tau)\to  (\Omega^\bullet(X),d)$ induces a quasi-isomorphism.
\ep

\begin{proof}[Proof of Corollary~\ref{cor}]
We observe that the values of the sheaves $(\Ell(\Map(T,X)),Q)$ and $(\Ell(X),d)$ on $U\subset \cH\to \Mell$ are chain complexes that compute the hypercohomology over~$U$: the higher derived functors are all trivial. By Theorem~\ref{thm:1}, this shows that the local values are the same, with restriction giving a quasi-isomorphism. This proves that $(\Ell(\Map(T,X)),Q)$ satisfies property (1). Property (2) is an immediate consequence of Theorem~\ref{thm:1} applied to $\cH/\SL_2(\Z)$. 
\ep


\section{Integration on the double loop space via localization}\label{sec:integral}


In this section, fix a Riemannian metric on~$X$. 

\subsection{An $\cH$-family of operators acting on the normal bundle to $T\subset \Map(T,X)$}

\begin{defn}\label{defn:norm} The (complexified) \emph{normal bundle} to the constant maps $X\subset \Map(T,X)$ is the infinite rank $C^\infty(X)$-module 
\beq
&&\Gamma(X;\nu):=\{v\in \Gamma(X\times T;p^*TX_\C)\mid \int_{T} v=0,\ p\colon X\times T\to X\}.\label{eq:Normal}
\eeq
\end{defn}
\begin{rmk}
One can also define $\Gamma(X;\nu)$ as the complexification of the sheaf of sections of $N\to X$ from Lemma~\ref{lem:openexists}, dropping the restriction that sections lie in an $\epsilon$-ball. 
\end{rmk}

\begin{defn} \label{defn:pairing}
Define the sheaf $\Ell(X;\nu):=\Ell(X)\otimes_{C^\infty(X)} \Gamma(X;\nu)$ on $\Mell$, viewing $\Gamma(X;\nu)$ as a constant sheaf. Hence, local values are
$$
\Ell(X;\nu)(U):=\Ell(X)(U)\otimes_{C^\infty(X)} \Gamma(X;\nu)\simeq \mathcal{O}(U;\Omega^\bullet(X;\C[\beta,\beta^{-1}])\otimes_{C^\infty(X)} \Gamma(X;\nu))
$$
for $U\subset \cH$, i.e., the constant sheaf valued in $\Omega^\bullet(X;\C[\beta,\beta^{-1}])\otimes_{C^\infty(X)} \Gamma(X;\nu)$. The $\SL_2(\Z)$-equivariant structure for $\Ell(X;\nu)$ is purely through the action on~$\beta$. 
\end{defn}

There is a $C^\infty(X)$-linear $T$-action on~$\Gamma(X;\nu)$ that gives operators $\mathcal{L}_{\partial_x}$ and $\mathcal{L}_{\partial_x}$ acting by derivations on $\Gamma(X;\nu)$ where $\partial_x,\partial_y$ are the standard basis for the Lie algebra of~$T$. We promote these to operators on sections of~$\Ell(X;\nu)$. Pulling back the curvature operator $R\in \Omega^2(X;\End(TX))$ along the projection $X\times T\to X$ defines a 2-form on~$X$ valued in endomorphisms of $\Gamma(X;\nu)$. This defines an endomorphism $\beta R\colon \Ell(X;\nu)\to \Ell(X;\nu)$ of sheaves. 

\begin{defn} 
The family of operators 
$$
\mathcal{L}_{\omega^{\sharp}}+\beta R:=\tau\mathcal{L}_{\partial_x}-\mathcal{L}_{\partial_y}+\beta R
$$
defines an endomorphism of the sheaf $\Ell(X;\nu)$ over $\Mell$. 
\end{defn} 


Finally, we note that there is a $T$-invariant $C^\infty(X)$-valued bilinear pairing on $\Gamma(X;\nu)$,
\beq
( v,w) =\int_{T}\langle v,w\rangle\in C^\infty(X),\quad v,w\in \Gamma(X;\nu)\label{eq:pairing}
\eeq
where $\langle-,-\rangle$ is the $\C$-bilinear extension of the Riemannian metric on $X$ pulled back along the projection $X\times T\to X$ and the integral uses the standard volume form on $T=\R^2/\Z^2$. The formula~\eqref{eq:pairing} extends to a pairing on $\Ell(X;\nu)$ with values in $\Ell(X)$.


%
%



\subsection{The regularized Pfaffian as a derived global section} 

Let $\Z^2_*=\Z^2\setminus (0,0)$. Each $(m,n)\in \Z^2_*$, defines an inclusion of $C^\infty(X)$-modules,
$$
f_{n,m}\colon \Gamma(TX_\C)\hookrightarrow \Gamma(X;\nu)\qquad v\mapsto e^{2\pi i (n x+my)}v
$$
for $e^{2\pi i (nx+my)}\in C^\infty(T)$. Endowing $TX_\C$ with the $T$-action of weight $(n,m)$ on its fibers, $f_{n,m}$ is $T$-equivariant, and hence the operators $\mathcal{L}_{\omega^{\sharp}}+\beta R$ restrict as
$$
f^*_{n,m}(\mathcal{L}_{\omega^{\sharp}}+\beta R)=2\pi i (n\tau-m) \id_{TX_\C}+\beta R\in \mathcal{O}(\cH;\Omega^\bullet(X;\End(TX_\C))).
$$
Furthermore, $\mathcal{L}_{\omega^{\sharp}}+\beta R$ is determined by these restrictions. 
The pairing in Definition~\ref{defn:pairing} also restricts along $f_{n,m}\oplus f_{-n,-m}\colon \Gamma(TX_\C)\oplus \Gamma(TX_\C)\to \Gamma(X;\nu_\C)$ and 
$$(\mathcal{L}_{\omega^{\sharp}}+\beta R)_{\pm m,\pm n}:=(f_{n,m}\oplus f_{-n,-m})^*(\mathcal{L}_{\omega^{\sharp}}+\beta R)\in \mathcal{O}(\cH;\Omega^\bullet(X;\End(TX_\C\oplus TX_\C)))$$
is skew symmetric.
Hence, $(\mathcal{L}_{\omega^{\sharp}}+\beta R)_{\pm m,\pm n}$ has a well-defined Pfaffian valued in $\mathcal{O}(\cH;\Omega^\bullet(X))$. 
We normalize with respect to the Pfaffian for the analogous operator with $R=0$, i.e., 
$$
\det{}^{1/2}_{\rm norm}\Big((\mathcal{L}_{\omega^{\sharp}}+\beta R)_{\pm m,\pm n}\Big):=\frac{\Pf((\mathcal{L}_{\omega^{\sharp}}+\beta R)_{\pm m,\pm n})}{\Pf(\mathcal{L}_{\omega^{\sharp}}|_{\pm m,\pm n})}\in \mathcal{O}(\cH;\Omega^\bullet(X)).
$$
We then define the regularized Pfaffian as the product of normalized Pfaffians. 

\begin{defn}
Let 
$$
\det{}^{\! 1/2}_{\rm ren}(\mathcal{L}_{\omega^{\sharp}}+\beta R):=\prod_{(n,m)\in \Z^2_+} \det{}^{\! 1/2}_{\rm norm}\Big((\mathcal{L}_{\omega^{\sharp}}+\beta R)_{\pm m,\pm n}\Big)
$$
where the product is over $\Z^2_+=\{(n,m)\in \Z^2\mid m<0 \ {\rm or} \ m=0, n>0\}.$ 
\end{defn}

\begin{rmk} One can also define the regularized determinant using methods of $\zeta$-regularization. \end{rmk} 

\begin{proof}[Proof of Theorem~\ref{thm:int}] 
For finite-dimensional Pfaffians, we have the formula 
$$\Pf\left(\begin{array}{cc} 0 & A \\ -A^T & 0\end{array}\right)=(-1)^{d(d-1)/2}\det(A).$$
We use this to compute
\beq
\det{}^{\!1/2}_{\rm norm}\Big((\mathcal{L}_{\omega^{\sharp}}+\beta R)_{\pm m,\pm n}\Big)&=&\frac{\det(2\pi i (n\tau-m) \id_{TX_\C}+\beta R)}{\det(2\pi i (n\tau-m))}=\det\Big(\id_{TX_\C}+\frac{\beta R}{2\pi i (n\tau-m)}\Big)\nonumber\\
&=&\exp\left(\sum_{k=1}^\infty \frac{\beta^{2k} {\rm Tr}( R^{2k})}{(2\pi i )^{2k} k(n\tau-m)^{2k}}\right).\nonumber
\eeq
The product over $(n,m)\in \Z^2_+$ is
\beq
\det{}^{\!1/2}_{\rm ren}(\mathcal{L}_{\omega^{\sharp}}+\beta R)&=&\prod_{(n,m)\in \Z^2_+} \exp\left(\sum_{k=1}^\infty \frac{\beta^{2k} {\rm Tr}(R^{2k})}{(2\pi i )^{2k} k(n\tau-m)^{2k}}\right)\nonumber \\
&=&\exp\left( \sum_{k=1}^\infty \sum_{(n,m)\in \Z^2_*} \frac{\beta^{2k} {\rm Tr}(R^{2k})E_{2k}}{(2\pi i )^{2k} 2k}\right)\label{eq:Witten}
\eeq
where $E_{2k}(\tau)=\sum\frac{1}{(n\tau+m)^{2k}}$ is the $2k$th Eisenstein series, which is absolutely convergent for $k>1$ and conditionally convergent for $k=1$. Fix a choice of ordering for this sum, and denote the result by~$E_2\in C^\infty(\cH)$. The result~\eqref{eq:Witten} is the function
$$
\Wit(X):=\exp\left( \sum_{k=1}^\infty \sum_{(n,m)\in \Z^2_*} \frac{\beta^{2k} {\rm Tr}(R^{2k})E_{2k}}{(2\pi i )^{2k} 2k}\right)\in C^\infty(\cH;\Omega^\bullet(X)[\beta,\beta^{-1}]).
$$
This fails to gives an $\SL_2(\Z)$-invariant holomorphic function precisely because there is no choice of ordering for which $E_2$ is a modular form. The failure of modularity and holomorphy are measured by the equivariant and Dolbeault differentials,
\beq
&&\delta\Wit(X) \in C^\infty(\cH\times \SL_2(\Z);\Omega^\bullet(X)[\beta,\beta^{-1}]),\quad \overline{\partial}\Wit(X) \in \Omega^{0,1}(\cH;\Omega^\bullet(X)[\beta,\beta^{-1}]).\label{eq:fail2}
\label{eq:fail1}
\eeq
where $\delta\Wit(X)=p^*\Wit(X)-{\rm act}^*\Wit(X)$ is the difference between the pullback along projection and action by $\SL_2(\Z)$. 
If $p_1(TX)=-{\rm Tr}(R^{2})/8\pi^2$ represents zero in cohomology, then $\delta\Wit(X)$ and $\overline{\partial}\Wit(X)$ are $d$-exact. A choice of $H\in \Omega^3(X)$ with $p_1(TX)=dH$ therefore determines a triple
$$
\mathbb{W}{\rm it}(X):=\left(\Wit(X),A_1,A_2\right)\qquad \delta\Wit(X)=dA_1, \quad \overline{\partial}\Wit(X)=dA_2
$$
which by~\eqref{eq:fail1} specifies a cocycle in the hypercohomology of~$(\Ell(X),d)$ over $\Mell$ of total degree zero, i.e., $\mathbb{W}{\rm it}(X)$ determines a derived global section. 
The class in $\TMF(X)\otimes \C$ underlying this derived global section is the Witten class of~$X$. If instead $p_1(TX)$ is not exact, then the class underlying $\Wit(X)$ does not define a class in cohomology with values in modular forms. Hence, by Proposition~\ref{prop:1} there is no way to extend $\Wit(X)$ to a derived global section of $\Ell(X)$. 
Finally, by the Riemann--Roch theorem, the pushforward 
$$
\langle f\rangle =\int_X f\cdot \Wit(X)\in \Gamma(\Mell,\Ell^{{\rm deg}(f)-{\rm dim}(X)}(\pt))
$$
is the one inherited from complexifying the string orientation of TMF. In particular $\langle 1\rangle=\Wit(X)$ is the Witten genus of $X$, a modular form of weight ${\rm dim}(X)/2$. 
\ep
\appendix

\section{The mapping space $\Map(T,X)$ as a sheaf}\label{sec:appen}

Let $T=\R^2/\Z^2$ be the standard 2-dimensional torus. 

\begin{defn}
Let $\Mfld$ be the category of smooth manifolds and smooth maps and $\TMfld$ the category of $T$-manifolds and $T$-equivariant maps. Let $\Sh(\Mfld)$ denote the category of (set-valued) sheaves on the site of manifolds with the usual open covers and $\Sh(\TMfld)$ the category of (set-valued) sheaves on $\TMfld$ with $T$-invariant open covers. 
\end{defn}

\begin{ex}
Any manifold (respectively, $T$-manifold) defines a representable sheaf on $\Mfld$ (respectively, $\TMfld$). We often identify a ($T$-)manifold with its associated sheaf. 
\end{ex}

For $\X\in \Sh(\Mfld)$ or $\Sh(\TMfld)$, let $\X(S)$ denote the value of $\X$ on~$S\in \Mfld$ or $S\in \TMfld$. By the Yoneda lemma, this set is in bijection with the maps of sheaves $S\to \X$, viewing $S$ as a representable sheaf. We also refer to $\X(S)$ as the \emph{$S$-points of $\X$}. 

\begin{defn}
For a graded commutative $\C$-algebra $A$, let $\Omega^\bullet(-;A)$ be the sheaf on $\Mfld$ that assigns to $S$ the algebra of $A$-valued differential forms,~$\Omega^\bullet(S;A)$; we use the same notation for this sheaf on the site of~$\TMfld$. Let $\Omega^\bullet(-;A)^T$ be the sheaf on $\TMfld$ that assigns the set of $T$-invariant differential forms, $\Omega^\bullet(S;A)^T$. Finally, set $\Omega^\bullet(-):=\Omega^\bullet(-;\C)$ and $\Omega^\bullet(-)^T:=\Omega^\bullet(-;\C)^T$.
\end{defn}

For $\xi$ a complex vector field on~$T$, there are morphisms in~$\Sh(\TMfld)$
\beq
&&d\colon \Omega^\bullet\to \Omega^{\bullet+1}\qquad \iota_\xi \colon \Omega^\bullet\to \Omega^{\bullet-1},\qquad \mathcal{L}_\xi\colon \Omega^\bullet\to \Omega^\bullet\label{eq:Cartan}
\eeq
from the de~Rham differential, contraction with the complex vector field associated with $\xi$ under the $T$-action, and the Lie derivative with respect to this vector field. These natural transformations satisfy the Cartan identity, $[d,\iota_\xi]=\mathcal{L}_\xi$. We also get a square zero natural transformation $d-\iota_\xi$ from $\Omega^\bullet(-)^T$ to itself which has odd parity when viewing $\Omega^\bullet(S)^T$ as a $\Z/2$-graded algebra. In $\Sh(\Mfld)$, we have the morphism $d\colon \Omega^\bullet\to \Omega^{\bullet+1}$. 

\begin{defn}
For a sheaf $\X$ on $\Mfld$ or $\TMfld$, \emph{$A$-valued differential forms on $\X$}, denoted $\Omega^\bullet(\X;A)$, are the set of maps $\X\to \Omega^\bullet(-;A)$. \emph{Functions on $\X$} and \emph{$T$-invariant forms on $\X$} are defined similarly and denoted by $C^\infty(\X)$ and $\Omega^\bullet(\X;A)^T$, respectively. 
\end{defn}

We observe that~$\Omega^\bullet(\X)=\prod_j \Omega^j(\X)$ and $\Omega^\bullet(\X)^T=\prod_j \Omega^j(\X)^T$ in our definition. 
%
For $\X\in \Sh(\Mfld)$, $(\Omega^\bullet(\X),d)$ is a cdga. For $\X\in \Sh(\TMfld)$ and $\xi$ a vector field on~$T$, $(\Omega^\bullet(\X;A)^T,d-\iota_\xi)$ is a $\Z/2$-graded complex with graded commutative multiplication. 

\begin{defn} For any smooth manifold $X$, let $\Map(T,X)$ denote the sheaf on $\Mfld$ whose $S$-points are maps $S\times T\to X$. We use the same notation for the sheaf on $\TMfld$ whose $S$-points are $T$-invariant maps $S\times T\to X$.
\end{defn}

\begin{rmk} The above definition is designed to satisfy the adjunctions, ${\rm Hom}(S,\Map(T,X))\simeq {\rm Hom}(S\times T,X)$ in $\Sh(\Mfld)$ or $\Sh(\TMfld)$ with the trivial $T$-action on $X$. \end{rmk}

\begin{defn}
A morphism $p\colon \mathcal{U}\to \X$ in $\Sh(\Mfld)$ or $\Sh(\TMfld)$ is \emph{open} if for all maps $S\to \X$ of sheaves where $S$ is representable, the pullback 
\beq
\begin{tikzpicture}[baseline=(basepoint)];
\node (A) at (0,0) {$U$};
\node (B) at (2,0) {$\mathcal{U}$};
\node (C) at (0,-1) {$S$};
\node (D) at (2,-1) {$\X$}; 
\draw[->] (A) to (B);
\draw[->] (A) to (C);
\draw[->] (C) to (D);
\draw[->] (B) to node [right] {$p$} (D);
\path (0,-.75) coordinate (basepoint);
\end{tikzpicture}\nonumber
\eeq
is a representable sheaf $U$ and the associated map $U\to S$ of manifolds is open. 
\end{defn}

\begin{defn} \label{def:homotopy}A \emph{homotopy} between maps of sheaves $f_0,f_1\colon \X\to \Y$ is a morphism $F\colon \X\times \R\to \Y$ whose restrictions to $\X\times \{0\}$ and $\X\times \{1\}$ are the maps~$f_0$ and~$f_1$, respectively, and where in $\Sh(\TMfld)$ we take the trivial $T$-action on $\R$. For a subsheaf~$\Y\subset \X$ a \emph{deformation retraction} is a homotopy $F\colon \X\times \R \to \X$ such that $f_0$ is the identity on $\X$, $f_1$ factors through~$\Y\subset \X$, and the restriction of $F$ to $\Y$ factors through the identity map on~$\Y$. 
\end{defn}

\bibliographystyle{amsalpha}
\bibliography{references}

\providecommand{\bysame}{\leavevmode\hbox to3em{\hrulefill}\thinspace}
\providecommand{\MR}{\relax\ifhmode\unskip\space\fi MR }
\providecommand{\MRhref}[2]{%
  \href{http://www.ams.org/mathscinet-getitem?mr=#1}{#2}
}
\providecommand{\href}[2]{#2}
\begin{thebibliography}{{Hai}08}

\bibitem[AB84]{AtiyahBott}
M.~Atiyah and R.~Bott, \emph{The moment map and equivariant cohomology},
  Topology \textbf{23} (1984).

\bibitem[AG83]{Alvarez}
L.~Alvarez-Gaum{\'e}, \emph{Supersymmetry and the {Atiyah-Singer} index
  theorem}, Communications in Mathematical Physics \textbf{90} (1983),
  161--173.

\bibitem[AHR10]{AHR}
M.~Ando, M.~Hopkins, and C.~Rezk, \emph{Multiplicative orientations of
  {KO}-theory and the spectrum of topological modular forms}, preprint (2010).

\bibitem[AHS01]{AHS}
M.~Ando, M.~Hopkins, and N.~Strickland, \emph{Elliptic spectra, the theorem of
  the cube, and the {Witten} genus}, Invent. math. \textbf{146} (2001),
  595--687.

\bibitem[AS02]{AlvarezSinger}
O.~Alvarez and I.~Singer, \emph{Beyond the elliptic genus}, Nuclear Physics B
  \textbf{633} (2002).

\bibitem[Ati71]{AtiyahSpin}
M.~Atiyah, \emph{Riemann surfaces and spin structures}, Annales scientifiques
  de l'\'Ecole Normale Sup\'erieure, S\'erie 4 \textbf{4} (1971).

\bibitem[Ati85]{AtiyahCircular}
\bysame, \emph{Circular symmetry and stationary phase approximation},
  Proceedings of the conference in honor of {L. Schwartz, Ast\'erisque}
  \textbf{131} (1985).

\bibitem[Ati88]{AtiyahTQFT}
\bysame, \emph{Topological quantum field theories}, Publications
  Math\'ematiques de {l'IH\'ES} \textbf{68} (1988).

\bibitem[Bis85]{Bismut3}
J.~Bismut, \emph{Index theorem and equivariant cohomology on the loop space},
  Comm. Math. Phys. \textbf{98} (1985).

\bibitem[Bis86]{Bismut1}
\bysame, \emph{Localization formulas, superconnections, and the index theorem
  for families}, Commun. Math. Phys. \textbf{103} (1986).

\bibitem[Bis11]{Bismut2}
\bysame, \emph{{Duistermaat-Heckman} formulas and index theory}, Geometric
  Aspects of Analysis and Mechanics \textbf{292} (2011).

\bibitem[BV81]{BV1}
I.~Batalin and G.~Vilkovisky, \emph{Gauge algebra and quantization}, Phys.
  Lett. B \textbf{102} (1981).

\bibitem[BV82]{BerlineVergneEquiv}
N.~Berline and M.~Vergne, \emph{Classes caracteristiques equivariantes.
  {Formule} de localisation en cohomologie equivariante}, Comptes rendus de
  l'Acad\'emie des sciences (1982).

\bibitem[BV83a]{BV2}
I.~Batalin and G.~Vilkovisky, \emph{Quantization of gauge theories with
  linearly dependent generators}, Phys. Rev. D \textbf{28} (1983).

\bibitem[BV83b]{BerlineVergne}
N.~Berline and M.~Vergne, \emph{Z\'eros d'un champ de vecteurs et classes
  caract\'eristiques \'equivariants}, Duke Math J. \textbf{50} (1983).

\bibitem[BV85]{BerlineVergne2}
\bysame, \emph{The equivariant index and {Kirillov} character formula}, Amer.
  J. of Math \textbf{107} (1985).

\bibitem[BX11]{BehrendXu}
K.~Behrend and P.~Xu, \emph{Differentiable stacks and gerbes}, J. Symplectic
  Geom. \textbf{9} (2011), no.~3.

\bibitem[Cra03]{Crainic}
M.~Crainic, \emph{Differentiable and algebroid cohomology, van {Est}
  isomorphisms, and characteristic classes}, Comment. Math. Helv. \textbf{78}
  (2003).

\bibitem[DH82]{DuistermaatHeckman}
J.~Duistermaat and G.~Heckman, \emph{On the variation in the cohomology of the
  symplectic form of the reduced phase space}, Inventiones Mathematicae
  \textbf{69} (1982).

\bibitem[GH04]{GoerssHopkins}
P.~Goerss and M.~Hopkins, \emph{Moduli spaces of commutative ring spectra},
  Structured ring spectra, London Math. Soc., vol. 315, Cambridge Univ. Press,
  Cambridge, 2004.

\bibitem[Goe09]{Goerss}
P.~Goerss, \emph{Topological modular forms (after {Hopkins, Miller, and
  Lurie})}, S\'eminaire {B}ourbaki (2009).

\bibitem[{Hai}08]{hain}
R.~{Hain}, \emph{{Lectures on Moduli Spaces of Elliptic Curves}}, ArXiv
  e-prints (2008).

\bibitem[Han08]{Han}
F.~Han, \emph{Supersymmetric {QFT}, super loop spaces and the {Bismut-Chern}
  character}, Preprint, 2008.

\bibitem[Hop94]{HopkinsICM94}
M.~Hopkins, \emph{Topological modular forms, the {Witten} genus, and the
  theorem of the cube}, Proceedings of the ICM \textbf{1} (1994).

\bibitem[JP90]{JonesPetrack}
J.~Jones and S.~Petrack, \emph{The fixed point theorem in equivariant
  cohomology}, Trans. Amer. Math. Soc. \textbf{322} (1990).

\bibitem[Lur09]{Lurie_Elliptic}
J.~Lurie, \emph{A survey of elliptic cohomology}, Algebraic Topology (N.~Baas,
  E.~Friedlander, J.~Bj\"orn, and P.~{\O}st\ae{r}, eds.), Abel Symposia,
  vol.~4, Springer Berlin Heidelberg, 2009 (English).

\bibitem[PMZ16]{localizationvol}
V.~Pestun and eds. M.~Zabzine, \emph{Localization techniques in quantum field
  theories}, Journal of Physics {A} (2016).

\bibitem[PZ16]{localizationvolintro}
V.~Pestun and M.~Zabzine, \emph{Introduction to localization in quantum field
  theory}, Journal of Physics {A} (2016).

\bibitem[Seg88]{Segal_Elliptic}
G.~Segal, \emph{Elliptic cohomology}, S\'eminaire N. Bourbaki \textbf{695}
  (1988).

\bibitem[Seg04]{Segal_CFT}
\bysame, \emph{The definition of conformal field theory}, Topology, Geometry,
  and Quantum Field Theory (2004).

\bibitem[ST04]{ST04}
S.~Stolz and P.~Teichner, \emph{What is an elliptic object?}, Topology,
  geometry and quantum field theory, London Math. Soc. LNS 308, Cambridge Univ.
  Press (2004), 247--343.

\bibitem[ST11]{ST11}
\bysame, \emph{Supersymmetric field theories and generalized cohomology},
  Mathematical Foundations of Quantum Field and Perturbative String Theory ({B.
  Jur{\v c}o, H. Sati, U. Schreiber}, ed.), Proceedings of Symposia in Pure
  Mathematics, 2011.

\bibitem[Wit82]{susymorse}
E.~Witten, \emph{Supersymmetry and {Morse} theory}, Journal of Differential
  Geometry \textbf{17} (1982), 661--692.

\bibitem[Wit88]{witten_dirac}
\bysame, \emph{The index of the {Dirac} operator in loop space}, Elliptic
  Curves and Modular Forms in Algebraic Topology (P.~Landweber, ed.), Springer
  Berlin Heidelberg, 1988.

\end{thebibliography}

\end{document}